\documentstyle[12pt]{article}
\begin{document}
\title{A structure and representations of diffeomorphism groups
of non-Archimedean manifolds.
\thanks{ Mathematics subject classification (1991 Revision):
43A65, 46S10, 57S05.} }
\author{ Sergey V. Ludkovsky}
\date{31 March 2000}
\maketitle
\par permanent address: 
Theoretical Department, Institute of General Physics, \\
Russian Academy of Sciences, \\
Str. Vavilov 38, Moscow, 117942, Russia.
\begin{abstract}
\par In this article diffeomorphism groups $G$
of manifolds $M$ on locally $\bf F$-convex spaces
over non-Archimedean fields $\bf F$ are investigated.
It is shown that their structure has many differences
with the diffeomorphism groups of real and complex manifolds.
It is proved that $G$ is not a Banach-Lie group,
but it has a neighbourhood $W$ of the unit element
$e$ such that each element $g$ in $W$ belongs to
at least one corresponding one-parameter subgroup.
\par It is proved that $G$ is simple and perfect.
Its compact subgroups $G_c$ are studied such that
a dimension over $\bf F$ of its tangent space 
$dim_{\bf F}T_eG_c$ in $e$ may be infinite.
This is used for decompositions of continuous 
representations into irreducible and investigations of
induced representations.
\end{abstract}
\section { Introduction.}
\par  This article is devoted to the investigation of a structure
and representations of diffeomorphism groups of non-Archimedean manifolds. 
In previous works \cite{lutmf99,lu4}
quasi-invariant measures on diffeomorphism groups relative to
dense subgroups were constructed. Irreducible representations
associated with the quasi-invaraint measures on groups and the corresponding 
configuration spaces were constructed in \cite{lutmf99,lusmb98}.
Classical diffeomorphism groups (that is, for real or complex manifolds) 
play very important role in hydrodynamics, quantum mechanics 
and superstring theory \cite{bao,ebi,isham}.
On the other hand, non-Archimedean quantum mechanics develops rapidly
\cite{jang,khr,vla3}. It is helpful in special situations, when
series or integrals divergent in quantum mechanics over the 
complex field $\bf C$ are convergent in the non-Archimedean case.
In particular, non-Archimedean diffeomorphism groups can be used
in non-Archimedean quantum mechanics and quantum gravity 
\cite{lutmf99,vla3}. 
\par There are many principal differences between classical and 
non-Archimedean functional analysis \cite{bosc,roo,sch1}.
This is the source why non-Archimedean diffeomorphism groups differ
in many respects from that of classical one.
\par In \cite{mat1} it was shown that classical diffeomorphism 
groups are simple and perfect, but proofs there are based
on local connectedness, homotopies,
the existence and the uniquiness of solutions of differential equations
in spaces of functions of the class of smoothness $C(t)$ 
for $t<\infty $. In the non-Archimedean case
even for the class of smoothness $C(\infty )$ there is not any uniquiness,
because of locally constant additional terms.
In the classical case the small inductive dimension
$ind (G)>0$ (for real manifolds $ind (G)=\infty $), but
in the non-Archimedean case $ind (G)=0$. Therefore, the proof 
of simplicity and perfectness 
in this paper differ principally from the classical case.
For compact complex manifolds the diffeomorphism groups are Lie groups
\cite{koba}, but in the non-Archimedean case, as it is proved below, 
it is untrue.
\par This article is devoted to more general diffeomorphism
groups than in \cite{lutmf99,lusmb98}.
Here are considered manifolds not only on Banach spaces
over local fields, but also on locally $\bf F$-convex spaces,
where $\bf F$ is an infinite field of characteristic $char ({\bf F})=0$
with non-trivial non-Archimedean valuation. Classes of smoothness $C(t)$ of 
manifolds $M$ considered below are $1\le t\le \infty $
and also analytic $t=an$ such that they are certainly not less than that
of $G$. In particular this encompasses the class of manifolds treated
by rigid analytic geometry (see about it in \cite{bosc,freput}).
This geometry is helpful in non-Archimedean superstring theory
and theory of homologies and cohomologies, but it is related with
very narrow class of analytic functions \cite{boum}. It is also 
extremely restrictive for non-Archimedean functional analysis and 
quantum theory. Therefore, differentiable manifolds of classes
$C(t)$ for $1\le t\le \infty $ also are considered below.
Historically spaces of classes $C(t)$ with $t\in \bf N$
had appeared in \cite{sch1,sch3} several 
years later after the use of analytic spaces and manifolds in
\cite{boum,tate}. Schikhof had used difference quotients 
of functions, Tate had used a topology stronger than 
the Zariski topology.
\par For locally compact groups there is a theory of induced 
representations from subgroups \cite{barut,fell}, but its 
development for non-locally compact groups meets serious problems,
because the case of non-locally compact groups is more complicated
\cite{lubp}.
In this article with the help of structural theorems of
diffeomorphism groups induced representations are investigated.
\par   In \S 2 definitions, notations and preliminary results are given.
In \S 3 the structure of diffeomorphism groups $Diff(t,M)$ is studied,
where $Diff(t,M):=Hom(M)\cap C(t,M\to M)$, 
$C(t,M\to N)$ is a manifold of $C(t)$-mappings
from a manifold $M$ into a manifold $N$ over the same field $\bf F$.
Besides classes $C(t)$ also classes $C_0(t)$
are considered over local fields $\bf K$. 
If $dim_{\bf F}M\ge \aleph _0,$ then $C(t,M\to M)$ is of non-separable 
type over $\bf F$, but $C_0(t,M\to M)$ is of separable type, when
${\bf F}=\bf K$ and $dim_{\bf K}M\le \aleph _0$.
Such groups $G(t,M):=Hom(M)\cap C_0(t,M\to M)$ are
helpful for the construction of quasi-invariant 
$\sigma $-finite measures.
The diffeomorphism group is investigated below as the topological 
group and as the manifold.
It is proved
that $Diff(t,M)$ are simple and perfect. Then its structure as a 
manifold is studied. 
Apart from manifolds $M$ on locally convex spaces $X$ 
over $\bf R$ \cite{pres} in the case of $X$ over $\bf F$ 
the existence of clopen (closed and open) subgroup 
$W$ in $Diff(t,M)$ is proved below such that for each
$g \in W$ there exists a one-parameter subgroup $<g^z:$ $z \in {\bf F}>$
to which $g$ belongs. Nevertheless, it is proved that $Diff(t,M)$
are not Banach-Lie groups. In \S 3 also families
of compact subgroups $ \{ G^n_{u,\bf K} \} $ of the group $G(t,M)$
are constructed such that
$\bigcup_{n,u,\bf K}G^n_{u,\bf K}$ is dense in $G(t,M)$. 
In the particular case of the local field ${\bf F}=\bf K$
such subgroups have the following property: the $\bf K$-linear span
$sp_{\bf K}(T_eG^n_{u,\bf K})$ of $T_eG^n_{u,\bf K}$ 
is dense in $T_eG(t,M)$.
This is the important difference from the case of $M$ 
on $X$ over $\bf R$ or $\bf C$, because the maximal compact subgroup
in $G(t,M)$ in the classical case may be only finite-dimensional
for finite-dimensional $X$ over $\bf R$ or $\bf C$ \cite{stra}.
This also is impossible in the classical case, when
$M$ is not a compact complex manifold. Embeddings of classical groups
into the diffeomorphism groups also are discussed,
because, for example, $Sp(2n,{\bf F})$ is very important for symplectic 
structures associated with Hamiltonians in quantum mechanics.
\par  In \S 4 continuous unitary representations and also representations 
in non-Archimedean Banach spaces are decomposed into irreducible.
Then induced representations are considered.
Moreover, two theorems (inductive-reductive and for
internal tensor product representations) about decompositions of
induced representations are proved. This opens new classes
of unitary representations.
\section{Topologies of non-Archimedean diffeomorphism groups.}
To avoid misunderstandings we first present our definitions and notations
in \S \S 2.1-2.4.
\par  {\bf 2.1. Remarks.} Let $\bf K$ be a local field,
that is, a finite algebraic extension of the field $\bf Q_p$ of
$p$-adic numbers and either $0\le t < \infty $ or $t=an$. 
Then $C_*(t,M\to N)$, $Diff(t,M),$ $G(t,M)$
and $GC(t,M)$ be the same spaces as in \cite{lutmf99,lusmb98}, 
where $M$ and $N$ are 
the corresponding Banach manifolds over $\bf K$, where
either $*=\emptyset $ ($\emptyset $ is omitted as the index)
or $*=0$ or $*=c$.
It is necessary to mention, that in in the case of $M$ with 
an infinite atlas spaces in \cite{lutmf99} are proper subspaces of
the corresponding spaces in \cite{lusmb98}. Then analogously we get
these spaces for the class of locally analytic functions with $t=la$.
Evidently, these spaces
are isomorphic for different choices of atlases
$At(M)$ and $At(N)$ for $M$ and $N$ of classes not less, 
than either $C(t)$ or $C_0(t)$ respectively, since 
the valuation group $\Gamma _{\bf K}:= \{ |x|: 0\ne x \in {\bf K} \} $
is discrete in $(0,\infty )$ and due to
\cite{luum985} and Lemma 7.3.6 \cite{eng} each atlas of $M$
or $N$ has a disjoint covering $At'(M)$ or $At'(N)$, 
which is a refinement of the initial covering. Indeed,
if $\phi : M\to M'$ and $\psi : N\to N'$ are 
$C_*(t)$-diffeomorphisms
(that is, $\phi $ is bijective and surjective
and $\phi \in C_*(t,M\to M')$, $\phi ^{-1} \in C_*(t,M'\to M)$,
analogously for $\psi $),
then $g\mapsto \psi \circ g \circ \phi ^{-1}$
is a diffeomorphism of $C_*(t,M\to N)$ with $C_*(t,M'\to N')$,
where $g\in C_*(t,M\to N)$.
\par {\bf 2.2. Notation.} Let $\bf F$ be an infinite 
field of characteristic $char ({\bf F})=0$
with a non-trivial non-Archimedean valuation.
For $b\in \bf R$, 
$0< b \le 1$, we consider the following mapping:
$$(1)\mbox{ }j_b(\zeta ):= \zeta ^b\in {\bf \Lambda _p} 
\mbox{ for }\zeta \ne 0,\mbox{ }  j_b(0):=0,$$ 
such that $j_b(*): {\bf F}\to \bf \Lambda _p$, where 
$\bf \Lambda _p$ is a spherically complete
field with a valuation group $ \{ |x|:$ $0\ne x \in {\bf \Lambda _p} \} =
(0,\infty )\subset \bf R$ such that ${\bf C_p}\cup {\bf F}
\subset \bf \Lambda _p,$ $\bf C_p$ denotes the field of complex numbers
with the valuation extending that of $\bf Q_p$ \cite{diar,roo,sch1,wei}. 
For a space $X$ with a metric $d$ in it let
$B(X,y,r):= \{ x\in X: d(x,y)\le r \} $
and $B(X,y,r^-):= \{ x\in X: d(x,y)< r \} $
denote balls in $X$, where $0<r.$
\par  {\bf 2.3. Definitions and Notes.} Let us consider locally convex 
spaces $X$ and $Y$ over $\bf F$. Suppose $F: U\to Y$ is a mapping, 
where $U\subset X$ is an open bounded subset.
The mapping $F$ is called differentiable if for each $\zeta \in \bf F$,
$x \in U$ and $h \in X$ with $x+\zeta h \in U$ there exists
a differential such that
$$(1)\mbox{ }DF(x,h):= 
dF(x+\zeta h)/d\zeta  \mid _{\zeta =0 } :=\lim_{\zeta \to 0} 
\{ F(x+\zeta h)-F(x) \}/\zeta $$  and
$DF(x,h)$ is linear by $h$, that is, $DF(x,h)=:F'(x)h$, where $F'(x)$
is a bounded linear operator (a derivative). Let 
$$(2)\mbox{ }\Phi ^b F(x;h;\zeta ):=(F(x+\zeta h)-F(x))/j_b(\zeta )\in 
Y_{\bf \Lambda _p}$$ be partial difference quotients of order $b$
for $0<b\le 1$, $x+\zeta h\in U$, $\zeta h\ne 0$, $\Phi ^0F:=F$, where
$Y_{\bf \Lambda _p}$ is a locally convex space obtained from $Y$ by 
extension of a scalar field from $\bf F$ to $\bf \Lambda _p$.
By induction using Formulas $(1-2)$
we define partial difference quotients of order $n+b$ for each $0<b\le 1$:
$$(3)\mbox{ }\Phi ^{n+b} F(x;h_1,...,
h_{n+1};\zeta _1,...,\zeta _{n+1} ):= 
\{ \Phi ^{n} F(x+\zeta _{n+1} h_{n+1};$$
$$h_1,..,h_n; \zeta _1,...,\zeta _n)- \Phi ^n F(x;h_1,...,h_n;
\zeta _1,...,\zeta _n) \}/j_b(\zeta _{n+1})$$ 
and derivatives
$F^{(n)}=(F^{(n-1)})'.$
Then $C(t,U\to Y)$ is a space of functions $F: U\to Y$ for which
there exist bounded continuous extensions
${\bar \Phi }^vF$ for each $x$ and $x+\zeta _ih_i \in U$
and each $0\le v \le t$, such that each derivative
$F^{(k)}(x): X^k\to Y$ is a continuous $k$-linear operator for each
$x\in U$ and $0<k\le [t],$ where $0\le t<\infty $,
$h_i\in V$ and $\zeta _i\in S:=B({\bf K},0,1),$
$[t]=n\le t$ and $\{ t\} =b$ are the integral 
and the fractional parts of $t=n+b$ respectively,
$U$ and $V$ are open neighbourhoods of $x$ and $0$ in $X$,
$U+V\subset U$.
In the locally $\bf F$-convex space $C(t,U\to Y)$ its uniformity is given by 
the following family of pseudoultranorms:
{\large
$$(4)\mbox{ }\| F\|_{C(t,U\to Y),u,w} 
:=sup_{(x, x+\zeta _ih_i\in U; h_i\in V; u(h_i)\ne 0; \zeta _i\in S;
i=1,...,s=[v]+sign \{ v\} ; 0\le v \le t)}$$ 
$$w \{ ({\bar \Phi }^v F)( x;h_1,..,h_s;\zeta _1,...,\zeta _s) \} /
[\prod_{i=1}^s u(h_i)]^v$$ }
where $0\le t\in \bf R,$ $sign (y)=-1$ for $y<0$, $sign (y)=0$ for $y=0$
and $sign (y)=1$ for $y>0$, $\{ u \} $ and $\{ w \} $
are families of pseudoultranorms in $X$ and $Y$ giving their 
ultrauniformities \cite{nari}.
\par Then the locally $\bf F$-convex space
$$(5)\mbox{ }C(\infty ,U\to Y):=\bigcap_{n=1}^{\infty }C(n, U\to Y)$$ 
is supplied with the ultrauniformity 
given by the family of pseudoultranorms $\| * \| _{C(n,U\to Y),u,w}$.
\par {\bf 2.4. Remarks.} Spaces of analytic functions 
$C(an_R, B(X,x,R)\to Y)$
of radius of converegence not less than $0<R$ are defined with the help
of convergent series of polylinear polyhomogeneous functions \cite{boum}
for normed spaces $X$ and $Y$ over $\bf F$. 
Spaces of locally analytic functions
$C(la,M\to Y)$ are defined as inductive limits of spaces 
$C(la_r,M\to Y)$ of locally analytic functions $f$ such that
for each $x\in M$ there exists its neighbourhood $U_x$ in $M$
for which $f|_{U_x}$ has an analytic extension on $B(X,x,r)$,
where $M\subset X$. Then using projective limits of normed spaces we can
construct $C(la,M\to Y)$ for locally $\bf F$-convex spaces $X$ and $Y$. 
\par For $C(m)$-manifolds $M$ and $N$ on locally $\bf F$-convex
spaces $X$ and $Y$ with atlases $At(M)=\{ (U_i,\phi _i):
i\in \Lambda _M\}$ and $At(N)=\{ (V_i,\psi _i): i\in \Lambda _N\}$
a mapping $F: M\to N$ is called of class $C(t)$ if $F_{i,j}$
are of class $C(t)$ for each $i$ and $j$, where $F_{i,j}=\psi _i\circ F
\circ \phi _j^{-1}$, $\infty \ge m\ge t\ge 0$,
$\phi _i: U_i\to \phi _i(U_i)\subset X$ and $\psi _l:
V_l\to \psi _l(V_l)\subset Y$ are diffeomorphisms,
$U_i,$ $V_l$, $\phi _i(U_i)$ and $\psi _l(V_l)$
are open in $M$, $N$, $X$ and $Y$ respectively,
$\phi _i\circ \phi _l^{-1}\in C(m,\phi _l(U_i\cap U_l)\to X)$ for
each $U_i\cap U_l\ne \emptyset $, analogously for $\psi _i$.
\par Let $\pi _{z_1,..,z_n} : X\to sp_{\bf F} \{ z_1,...,z_n \} $
be a projection, where $z_1,...,z_n$ are linearly independent
vectors in $X$, then we set $C_0(t,M\to N)$ to be a 
completion of a subspace of cylindrical functions $f$
of class $C(t)$, that is, for each such $f$ there are 
$n\in \bf N$ and $z_1,...,z_n$ linearly independent 
in $X$ and $h\in C(t,(M\cap sp_{\bf F} \{ z_1,...,z_n \} ) \to N)$
such that $f(x)=h(\pi _{z_1,...,z_n}(x))$.
If $\theta : M\to N$ is a fixed mapping,
then $C^{\theta }_*(M,Y)$ is a space of functions $f: M\to Y$
such that $(f-\theta ) \in C_*(t,M\to Y)$, that induces a space
$C^{\theta } _*(t,M\to N)$, where $*=\emptyset $ or $\theta =0$.
\par Certainly we suppose throughout the paper, 
that $M$ and $N$ are of class $C(\tau )$ 
for spaces $C_*(t,M\to N)$ such that
$\tau =\infty $ for $0\le t\le \infty $, $\tau =an_r$ for $t=an_r$,
$\tau =la$ for $t=la$. 
\par Then $Diff(t,M):=Hom(M)\cap C^{id}(t,M\to M)$
and $G(t,M):=Hom(M)\cap C^{id}_0(t,M \to M)$
denote diffeomorphism groups for $t\ge 1$ or $t=an_R$ or $t=la$
and a homeomorphism group for $0\le t <1$ analogosuly to
\cite{lusmb98}, where $Hom(M)$ is the standard homeomorphism
group of $C(0)$ bijective surjective mappings of $M$ onto itself,
where the manifold $M$ is on the locally $\bf F$-convex space $X$
for $t\ne an_R$ and $X$ is the normed space for $t=an_R.$ 
\par {\bf 2.5.} Let $H$ be a locally $\bf F$-convex space, 
where $\bf F$ is a non-Archimedean field. Let $M$ be a 
topological manifold modelled 
on $H$ and $At(M) = \{ (U_j,f_j):$ $ j\in A \} $ be an atlas
of $M$ such that $card (A)\le w(H)$, where $f_j: U_j\to V_j$
are homeomorphisms, $U_j$ are open in $M$, $V_j$ are open in $H$,
$\bigcup_{j\in A}U_j=M$, $f_i\circ f_j^{-1}$ are continuous on 
$f_j(U_i\cap U_j)$ for each $U_i\cap U_j\ne \emptyset $.
Let $\bf \tilde F$, $\tilde H$ and $\tilde M$ denote completions
of $\bf F$, $H$ and $M$ relative to their uniformities.
\par {\bf Theorem.} {\it If either $H$ is infinite-dimensional over $\bf F$,
or $\bf \tilde F$ is not locally compact, 
then $M$ is homeomorphic to the clopen subset of $H$.}
\par {\bf Proof.} Since $\tilde H$ is the complete
locally $\bf \tilde F$-convex space, then
${\tilde H}=pr-\lim \{ {\tilde H}_q, \pi ^q_v, \Upsilon \} $ 
is a projective limit
of Banach spaces ${\tilde H}_q$ over $\bf \tilde F$, where $q\in \Upsilon $, 
$\Upsilon $ is an ordered
set, $\pi ^q_v: {\tilde H}_q\to {\tilde H}_v$ 
are linear continuous epimorphisms.
Therefore, each clopen subset $W$ in $\tilde H$ has a decomposition
$W=\lim \{ W_q, \pi ^q_v, \Upsilon \} $, 
where $W_q=\pi ^q_v(W)$ are clopen in
${\tilde H}_q$. The base of topology of $\tilde M$ 
consists of clopen subsets.
If $W\subset V_j$, then $f_j^{-1}(W)$ has an analogous decomposition.
From this and Proposition 2.5.6 \cite{eng} it follows, that
$\tilde M=\lim \{ {\tilde M}_q, {\tilde \pi }^q_v, \Upsilon \} $, where 
${\tilde M}_q$ are manifolds
on ${\tilde H}_q$ with continuous bonding mappings between charts
of their atlases.
If $H$ is infinite-dimensional over $\bf F$, then each ${\tilde H}_q$ is 
infinite-dimensional over $\bf \tilde F$ \cite{nari}.
From $card(A)\le w(H)$ it follows, that each ${\tilde M}_q$ has an atlas
$At'({\tilde M}_q)=\{ {U'}_{j,q};f_{j,q};{A'}_q \} $ 
equivalent to $At({\tilde M}_q)$ such that
$card ({A'}_q)\le w(H_q)=w({\tilde H}_q)$, since $w({\tilde H})=w(H)$,
where $At({\tilde M}_q)$
is induced by $At({\tilde M})$ by the quotient mapping
${\tilde \pi }_q: {\tilde M}\to {\tilde M}_q$.
In view of Theorem 2 \cite{luum985} each ${\tilde M}_q$ is homeomorphic 
to a clopen subset ${\tilde S}_q$ of ${\tilde H}_q$, 
where $h_q: {\tilde M_q}\to {\tilde S_q}$ are homeomorphisms. 
To each clopen ball $\tilde B$ in ${\tilde H}_q$ there corresponds
a clopen ball $B={\tilde B}\cap H_q$ in $H_q$, hence
$S_q={\tilde S}_q\cap H_q$ is clopen in $H_q$ and
$h_q: M_q\to S_q$ is a homeomorphism.
Therefore, $M$ is homeomorphic to
a closed subset $V$ of $H$, where $h: M\to V$ is a homeomorphism,
$V\subset H$, $h=\lim \{ id, h_q, \Upsilon \} $, 
$id: \Upsilon \to \Upsilon $
is the identity mapping. Since each $h_q$ is surjective, then
$h$ is surjective by Lemma 2.5.9 \cite{eng}.
If $x\in M$, then ${\tilde \pi }_q(x)=x_q\in M_q$,
where $\pi : H\to H_q$ are linear quotient mappings
and ${\tilde \pi }_q: M\to M_q$ are induced quotient mappings.
Therefore, each $x\in M$ has a neighbourhood ${\tilde \pi }_q^{-1}(
Y_q)$, where $Y_q$ is an open neighbourhood of $x_q$ in $M_q$.
Therefore, $h(M)=V$ is open in $H$.
\par {\bf 2.6. Theorems.} {\it 
\par $\bf 1.$ The spaces $Diff(t,M)$, $G(t,M)$ and $GC(t,M)$ are the
topological groups. 
\par $\bf 2.$ They have embeddings as clopen subsets
into the spaces $C_*(t,M\to X)$, where either $*=\emptyset $
or $*=0$ or $*=c$ respectively.
\par $\bf 3.$ If $\bf F$ and $X$ are complete, then $Diff(t,M),$ 
$G(t,M)$ and $GC(t,M)$ are complete.
\par $\bf 4.$ $G(t,M)$ and $GC(t,M)$ 
are separable for separable $M$.
\par $\bf 5.$ $Diff(t,M)$, $G(t,M)$ and $GC(t,M)$ are ultrametrizable 
for a manifold $M$ with a finite atlas $At(M)$ 
on a normed space $X$ and either $0\le t<\infty $ or $t=an_r$.}
\par  {\bf Proof.} $(A).$ Using the projective limits of normed spaces we
can reduce the proof to the case of $M$ on a normed space $X$,
since for each continuous either linear mapping $A: X\to X$ or 
polylinear and polyhomogeneous mapping on $X$ 
there are a pseudoultranorm $u$ in $X$ and a continuous mapping
either linear $\mbox{}_uA(x+ker (u))=A(x)$
or polylinear and polyhomogeneous 
$\mbox{}_uA(x_1+ker (u),...,x_n+ker (u))=A(x_1,...,x_n)$ 
from $X_u$ into $X_u$, where 
$X_u:=X/ker (u)$, $x, x_1,...,x_n \in X$, $x+\ker (u) \in X_u$
(see Theorem (5.6.3) \cite{nari}). 
The second statement is the consequence of Theorem 2.5.
If $f, g \in Diff(t,M)$ such that $0<t$, then for each
$0<b\le \min (1,t)$ we have
$(\Phi ^bf\circ g)(x;\xi ;h)=(\Phi ^b f)(g(x);\zeta ;z)$,
where either $\zeta =\xi $ and $z=(\Phi ^1g)(x;\xi ;h)$
for $b=1$, or $\zeta \in \bf F$ and $z\in X$ such that
$\zeta z=g(x+\xi h)-g(x)$ and $|\xi |^b/p\le |\zeta |\le |\xi |$,
$p$ is a prime number such that ${\bf Q_p}\subset \bf F$.
In view of recurrence Relations $2.3.(3)$ we get that $Diff(t,M)$
is the topological group for each $0\le t\le \infty $.
In view of definitions $Diff(an_R,M)$ and $Diff(la,M)$ 
are also topological groups.
\par $(B).$ Let at first $At(M)$ be finite. If 
$(f_n: n)$ is a Cauchy net in $C_*(t,M\to Y)$,
then $(\Phi ^vf_n: n)$ are uniformly convergent sequences
for each $0\le v \le t$ and $0\le t\le \infty $, also for each $v$ while
$t=an_r$. Consequently, $\lim_{n\to
\infty } (\Phi ^vf_n)=:F^v$ $\in C_*(\tau ,M^{s+1}\to Y)$, 
where $\tau =0$ for $0\le t\le \infty $ or
$\tau =an_r$ for $t=an_r$, $s:=[v]+sign ( \{ v \} )$. 
\par The statement about ultrametrizability follows from
\S 2.4 \cite{lusmb98} and \S 2.2 \cite{lutmf99}.
If $X$ is the Banach space, then from the completeness of $C_*(t,M\to Y)$,
in which either $Diff(t,M)$ or $G(t,M)$ or $GC(t,M)$ respectively
are closed, it follows that the latter spaces are also complete
(see Theorems 8.3.6 and 8.3.20 \cite{eng}).
\par  $(C).$ In the case $G(t,M)\ni f,g$ for $0\le t \le \infty $
due to \S \S 2.1-2.4 \cite{lusmb98} 
there is the equality
$$f_{i,j}\circ g_{j,l}(x)=\sum_{i \in I, n \in I, m \in {\bf N_o^n} }
a(m,f^k_{i,j}) {\bar Q}_m((g_{j,l})_n(x))q_i ,$$ 
where $(g)_n=(g^{i(1)},..., g^{i(s)})$, $g_{j,l}=( g^k_{j,l}(x):$
$U_l\to {\bf K}| k \in I)$, $M$ is modelled on $X=c_0(I,{\bf K})$, the set
$ \{ i\in I:$ $m(i)\ne 0 \} =$ $ \{ i(1),..., i(s) \} $ is finite,
$f_{i,j}= \phi _i\circ f\circ \phi _j^{-1}$ with the
corresponding domains, $s \in \bf N$, $n=Ord(m)$, 
$${\bar Q_m}((g)_n)=\prod_{j=1}^sQ_{m(i(j))}(g^{i(j)}) \mbox{ and } 
Q_{m(i)}(g^i):=P_{m(i)}(g^i)/ P_{m(i)}(u(m(i)),$$
where $P_{m(i)}$ are polynomials.
\par Coefficients $a(m,f_{i,j}^k)={\tilde
\Delta }^m(f_{i,j}^k(x))|_{x=0}$ are given by Corollary 2 from Proposition
7 \cite{ami}.
The polynomials $[\bar Q_m(x):$ $|m|\le n, m(j)\ne 0\mbox{ for }
j \in (i(1),...,i(s))]$ may be expressed throughout
$[x^m:$ $|m|\le n, m(j)\ne 0,\mbox{ for }j\in (i(1),...,i(s))]$
and vice versa, where $x^m=\prod_{m(j)\ne 0}x(j)^{m(j)}$,
$x(j) \in \bf K$, $x \in
B(X,0,1)$. Therefore, $\tilde \Delta ^mS_l(x)|_{x=0}=0$ for each
polynomial $S_l(x)$ with $l=(l(i):$ $i\in I, {\bf N_o}\ni l(i)\le m(i))$.
Whence the coefficients $a(m,f^k_{i,j}\circ g_{j,l'})$ may be expressed
throughout $a(l,f_{i,j}^k)a(q_{i(1)},
g_{j,l'}^{i(1)})$ $...a(q_{i(s)},g_{j,l'}^{i(s)})R_{l,i,q}/P_l(\tilde 
u(l)),$ where
$$(ii)\mbox{ }k+|l|+Ord(l)+\sum_{j=1}^s(|q_{i(j)}|+Ord(q_{i(j)}))-s
\ge k+|m|+Ord(m),$$  $q=(q_{(i(j))}\in {\bf N_o^{Ord(q_{i(j)})}}:$
$j=1,...,s)$, $0\le s\le |l|$, $R_{l,i,q}$ are polynomials by $u(i',j')$,
that appear from the decomposition of $\bar Q_m((g_{j,l})_n)$ in the form
of sums of products of $(g_{j,l})^k$ and $u(i',j')$ divided by 
$\bar P_m(\tilde u(m))$.
In view of $(i, ii)$ we get that $f\circ g\in G(t,M)$
and continuity of the composition, since in $(ii)$ for $|m|+Ord(m)
\to \infty $ or $|l|+Ord(l)\to \infty $ or there is $q_{i(j)}$
with $|q_{i(j)}|+Ord(q_{i(j)})\ge [|m|+Ord(m)+1]/s$. At the same time
$s>0$ for large $|m|+Ord(m)$. For $f=g^{-1}$ we get recurrence relations
for $a(m,(f_{i,j}^{-1})^k)$ throughout $a(m,f_{i,j}^l)$. From them
follows that $\rho _0^t(f^{-1},id)$ are polynomials of the Bell type
by $\rho _0^{\kappa }(f,id)$ in $1/p$ neighbourhood of $id$, where
$\kappa =0,1,...,[t], t<\infty ,$
$\rho _0^t$ is an ultrametric in $G(t,M)$ (see also \cite{lusmb98} and
Chapter 5 \cite{rio}).
This gives $f^{-1}\in G(t,M)$ and continuity of
the inversion $f\to f^{-1}$.
The case $t=\infty $ follows from Formula $2.3.(5)$.
\par  $(D).$ Now let $t=an_r$ and using the transformation
$x\to x\xi $ with $|\xi |=1/r$ we restrict the consideration to $r=1$.
If $g \in Diff(an_1,M)$ (or $G(an_1,M)$),
then $\| g\| \le 1.$ Indeed, there are the natural embeddings $\theta :
B({\bf K^n},0,1)\to B(X,0,1)$, consequently,
there are the restrictions
$g|_{M_n}:=g(\theta (x_n))$, where $\theta (x_n)=(x \in B(X,0,1):$ $
\theta (x_n)(i)=x(j)\mbox{ for }i=i(j)\in (i(1),...,i(n)),
\theta (x_n)(i)=0\mbox{ in others cases })$, $M_n=M\cap
\theta (B({\bf K^n},0,1)$).
In view of \S 54.4 in \cite{sch1} with the help of \cite{ami}
we get that if
$$f(x)=\sum_{m \in {\bf No^n} } a(m,f){\bar Q}_m(x)
\in C(0,B({\bf K^n},0,1)\to {\bf K}),$$ 
then $f$ is analytic if and only if there exists 
$$(iii)\mbox{ }\lim_{\mid m \mid \to \infty} a(m,f)/P_m(\tilde u(m))=0.$$
Moreover, in $C(an_1,B({\bf K^n},0,1)\to {\bf K})$
the following norms
$$(iv)\mbox{ } \Vert f\Vert :=\sup \{ \mid a(m,f)\mid J(an,m):
m \in {\bf No^n} \}
\mbox{  }  \mbox{ and }$$ 
$$(v)\mbox{ } \Vert f\Vert ":=\sup\{\mid b(m,f)\mid : m \in {\bf No^n} \} $$
are equivalent, where $J(an,m):=\mid 1/P_m(\tilde u(m))\mid $,
$b(m,f)$ are expansion coefficients by $x^m$.
Each function $g^k(\theta (x_n))$
is analytic and depends from a finite number of variables.
If $\| g\| >1$, so there is $M_n$ with $\| g(\theta (x_n))\| >1$.
\par  The basis ${\bar Q}_m(x)$ is orthogonal in the non-Archimedean
sense on $B({\bf K^n},0,1)$ with $\| \bar Q_m\| _{C(0,B\to {\bf K})}=1$.
Hence $|g^k(\theta (x_n))|>1$
contradicts $g \in Hom(M)$ and $M\subset B(X,0,1)$.
Therefore, $|a(m,g^k(\theta (x_n)))|J(an,m)\le 1$
for each $k, n$ and such $m$, $\theta $.
Hence $\| g\| _{C_*(an_1,M\to M)}\le 1$, 
since $\theta $ has the natural extension $\theta :
{\bf K^n}\hookrightarrow X$ such that $\theta $ is linear on
$\bf K^n$ and it is the embedding. 
Therefore, the composition and the inversion operations 
are correctly defined and they are continuous 
in $Diff(an_1,M)$ and $G(an_1,M)$ due to Formulas $(i, ii)$.
\par $(E)$. Now let $At(M)$ be infinite.
If $\bf F$, $X$ and $Y$ are complete, 
then $C_*(t,M\to Y)$ is complete 
(due to theorem about strict inductive limits
in Chapter 12 \cite{nari})
for $t\ne la$. If $(f_{\gamma }: \gamma \in \alpha )$ is a
Cauchy net in $C_*(la, M\to Y)$, consequently, there exist
$\delta \in \alpha $, $E\in \Sigma $ and $r_0>0$ such that
$supp(f_{\gamma }) \subset U^E$ and $f_{\gamma }\in
C_*(an_{r_0},M\to Y)$ for each $\gamma >\delta $, since
$\Pi ^r_R: C_*(an_R,U^E\to {\bf K})$ $\hookrightarrow
C_*(an_r,U^E\to {\bf K})$ are compact operators
for each $0<r<R$, where $\alpha $ is a limit ordinal,
$\Sigma $ is a family of all finite subsets of $\Lambda _M$. From
the completeness of $C_*(an_{r_0},M\to Y)$ it follows that
$( f_{\gamma })$ converges in $C_*(la,M\to Y)$, hence $C_*
(la,M\to Y)$ is complete. From definitions it follows
that $G(t,M)$, $Diff(t,M)$ and $GC(t,M)$ are
closed in $C_*(t,M\to M)$ for $*=0$, $*=\emptyset $
and $*=c$ respectively, whence they are also complete.
\par $(F)$. For separable $M$ and
$N$ the spaces $C_*(t,U^E\to N)$ are separable for each $E\in \Sigma $.  
The space $C_*(t,M\to Y)$ is isomorphic with the quotient space
$Z/P$, where $Z=\bigoplus_{j\in \Lambda }C_*(t,U_j\to Y)$,
$P$ is closed and $\bf K$-linear in $Z$. From the separability of
$Z$ and $\Lambda \subset \bf N$ it follows that
$C_*(t,M\to Y)$ is separable.
\par  $(G)$. From formulas (i,ii) it follows
that $GC(t,M)$ is the topological group for $M$ with
the finite atlas. For $f$ and
$g\in C_*(t,M\to M)
\cap Hom(M)$ for $0\le t\le \infty $ or $t=an_r$ there
are $E(f)$ and $E(g)\in \Sigma $ for which $supp(f):=cl
\{ x\in M:$ $f(x)\ne x \} \subset U^{E(f)}$ and $supp(g)\subset
U^{E(g)}$. Considering $f(supp(f))$
and $g(supp(g))\subset M$ homeomorphic with $supp(f)$ and $supp(g)$
correspondingly we get $g^{-1}\circ f\in C_*(t,
M\to M)\cap Hom(M)$. If $(f_{\gamma }: \gamma \in \alpha )$
and $(g_{\gamma }: \gamma \in \alpha )$ are two convergent nets
in either $G(t,M)$ or $Diff(t,M)$ or $GC(t,M)$) to
$f$ and $g$ respectively, so for each neighbourhood $W\ni id$
there exist $E\in \Sigma $ and $\beta \in \alpha $ such that
$g_{\gamma }^{-1}\circ f_{\gamma } \in W\cap 
C_*(t,U^E\to M)\cap Hom(M)$ for $0\le t\le \infty $ or
$t=an_r$, where $\alpha $ is a limit ordinal. Therefore,
for such $t$ the mapping $(f,g)\to g^{-1}\circ f$ is
continuous in $G(t,M)$ or $Diff(t,M)$ or 
$GC(t,M)$ respectively.
\par  For $t=la$ let $r=\min (r(f), r(g))$, where $f$ and $g
\in C_*(la,M\to M)\cap Hom(M)$, that is, there exist
$r(f)$ and $r(g)\in \Gamma _{\bf F}$ such that
$f\in C_*(an_{r(f)},M\to M)\cap Hom(M)$ and analogously for $g$.
Then $r\in \Gamma _{\bf F}$ and $g^{-1}\circ f
\in C_*(an_r,M\to M)\cap Hom(M)$. If $(f_{\gamma }:
\gamma \in \alpha )$ converges to $f$ and $(g_{\gamma })$
to $g$, then for each
neighbourhood $W\ni id  $ in $C_*(la,M\to M)\cap Hom(M)$
there exist $\beta \in \alpha $ and $E\in \Sigma $ such that
$(supp(g_{\gamma }^{-1}\circ f_{\gamma }))\cup (supp(g_{\gamma }))
\cup (supp(f_{\gamma }))\subset U^E$ for each $\gamma >\beta $
and $r(g_{\gamma })\ge r$, $r(f_{\gamma })\ge r$. 
Therefore, $(g_{\gamma }\circ
f_{\gamma }:$ $\gamma \in \alpha )$ converges to $g^{-1}\circ f$
in $C_*(la,M\to M)\cap Hom(M)$, consequently, the last
space is the topological group.
\section { A structure of diffeomorphism groups.}
\par  {\bf 3.1. Theorem.} {\it Let the groups $G=Diff(t,M)$ and $G=G(t,M)$ 
be the same as in \S 2.4, where either $1\le t \le \infty $
or $t=an_r$ or $t=la$.  
\par $\bf (1).$ If $M$ is on a complete space $X$, then
there exists a clopen subgroup $W$ in $G$ such that,
each element $g \in W$ belongs to the corresponding
one-parameter subgroup. 
\par $\bf (2).$ $Diff(t,M),$ $G(t,M)$ and $GC(t,M)$ are not 
Banach-Lie groups.}
\par  {\bf Proof.} As in the proof of Theorem 2.6 we can
use the projective limit $X=pr-\lim X_u$ of normed spaces 
$X_u$ that reduce the proof to the case of the 
manifold $M$ on the normed space $X$. 
\par $\bf (1).$ Let at first $G=G(t,M)$ 
and $M$ be with a finite atlas on $X$ over a local field $\bf K$.
We put $W:=\{ f\in G: \rho ^{\tau }_0(f,id)
\le p^{-2} \}$, then each
$f \in W$ is an isometry of $M$, where $\tau =t$ for either 
$1\le t\ne \infty $ or $t=an_r$,
$\tau \in \bf N$ for $t=\infty $. If $f, g\in W$, then $\rho _0^{
\tau }(f\circ g,id)=\rho ^{\tau }_0(g,f^{-1})\le \max (\rho _0^{
\tau }(g,id), \rho _0^{\tau }(id, f^{-1}))=$ $\max (\rho _0^{\tau }
(g,id), \rho _0^{\tau }(f,id))$. Therefore, $W$ is the
subgroup in $G$.
\par  Let at first $M_n$ be finite-dimensional over $\bf K$. 
There exists a restriction $f|_{M_n}$ for each $f\in G$, 
where $M_n$ is an analytic submanifold,
$\theta  : M_n\hookrightarrow M$ is an embedding, 
$dim_{\bf K}M_n=n \in \bf N$
is a dimension of $M_n$ over $\bf K$.
Since, locally polynomial functions $f(x)=id(x)+P(x)$ are dense in
$W$, it is sufficient to prove that each such $f(x)$ belongs
to a one-parameter subgroup. Here
$deg P=m \in {\bf N}$ is a degree of a polynomial,
$x \in M$ are a local coordinates. Denote $f_{i,j}=\phi _i\circ f
\circ \phi_j^{-1}$ simply by $f$ and $U_j$ by $M$. Let
$$g(j;x)=\sum_{s=0}^{\infty } A(j;x)^s x/s!,\mbox{ where}$$
$$A(j;x):=\sum_{i=1}^n T(j,i;x)\partial_i ,$$
$T(j,i;x)$ are polynomials on the
$j$-th iteration, $A(j;x)^sx:=A(j;x)(A(j;x)^{s-1}x)$ for $s>1$,
$A(j;x)^0:=x$, $\partial _i:=\partial /\partial x^i$.
Suppose $T(0,i;x)=P^i(x)$ for $i=1,...,n$,
and $A(0;x)A(1;x)+A(1;x)={\bar P}(x)$, where 
$$P(x)=\sum_{i=1}^n P^i(x)e_i \mbox{ and }
{\bar P}(x)=\sum_{i=1}^n P^i(x)\partial_i.$$
For the coefficients $T(1,i;x)$ there is the system of linear algebraic
equations, that gives the unique solution $A(1;x)$ with
$$\| A(0;x)T(1,i;x)\| _{\tau }\le \| T(1,i;x)\| _{\tau } \times
\| A(0;x)\|_{\tau },$$ since
$\| A(0;x)\|\le \| P(x)\|_{\tau },$ where 
$$\| A(j;x)\|_{\tau }:=\sup_{g\ne 0, g\in C_0(\tau ,M\to X)}
\| A(j;x)g \|_{\tau '}/\| g\|_{\tau },$$ 
$\tau '=\tau -1$ for $1\le \tau
<\infty $, $\tau '=\tau $ for $\tau =an_r$, 
$$\| g\|_{\tau }:=\| g\|_{C_0(\tau ,M\to X)}.$$
\par  Therefore, 
$$\| P(x)\|_{\tau } =\max_{i=1,...,n} \| T(1,i;x)\|_{\tau } ,$$
since $\| A(0;x)\|_{\tau } \le p^{-2}$.
Moreover, 
$$\max_{i=1,...,n} \| T(0,i;x)-T(1,i;x)\|_{\tau } 
\le \| P\|_{\tau } /p^2 ,$$
since $\| A(0;x)T(1,i;x)\|_{\tau } \le \| T(1,i;x)\|_{\tau } /p^2$
for each $i$.
\par Further by induction let for $j>0$ are satisfied the
following conditions:
$$(i)\mbox{ } {\bar P}(x)=A(j;x)+ \sum_{s=1}^j A(j-1;x)^s A(j;x)/s! ,$$
$$(ii)\mbox{ } \max_{i=1,...,n} \| T(j,i;x)-T(j-1,i;x)\|_{\tau }
\le p^{-j} \| P\|_{\tau } \mbox{ and }$$
$$(iii)\mbox{ } max_{i=1,...,n} \| T(j,i;x) \|_{\tau }=
\| P\|_{\tau }.$$
For $j+1$ instead of $j$ there exists the unique solution
$A(j+1;x)$ of the equation
$(i)$, since $\bar P(x)=(I+S_{j+1})A(j+1,x)$ with $\| S_{j+1}\|
\le 1/p$, $I$ is the identity operator. To Equation $(i)$ there corresponds
the linear algebraic equation $(I_z+F)Z=Y$, $Z$ and $Y\in \bf K^z$,
$z \in \bf N$, $I_z$ is the unit matrix
and $F$ is a matrix of size $z\times z$, $F=(F_{i,j})_{i,j=1,...,z}$,
$F_{i,j}\in \bf K$, $\max_{i,j}|F_{i,j}|\le 1/p$, $|det(I_z+F)|=1$.
Then 
$$\| A(j;x)^s T(j+1,i;x)/s!\| \le  \| T(j+1,i;x) \| p^{-s(2-1/(p-1))}
\mbox{ and}$$ 
$$t':=\max_{i=1,...,n} \| T(j+1,i;x)\|=\| P\| ,$$ since $\| A(j;x)\|
\le p^{-2}$. Consequently, 
$$\| [A(j+1;x)^s -A(j;x)^s]/s!\| \le
\| A(j+1;x)-A(j;x)\| p^z,$$
where $z=-(s-1)(2-1/(p-1))$, since 
$[A^i,B]=\sum_{l=0}^{i-1} A^l[A,B]A^{i-1-l},$ \\
$\| [A,B]\| \le \max \{ \| AB\| , \| BA\| \} $, $[A,B]:=AB-BA$, \\
$A^i-B^i=A(A^{i-1}+A^{i-2}B+...+B^{i-1})
-(A^{i-1}+A^{i-2}B+..+B^{i-1})B,$ \\
$\| AB\| \le \| A\| \times \| B\| $. Taking $A=A(j+1;x)$, $B=A(j;x)$
and using Formulas $(ii, iii)$ we get $\| AB-BA\| \le \max \{ \| AB-B^2 \|,
\| B^2-BA\| \}$ $\le \| A-B\|
/p^2$. From this it follows that 
$$\max_{i=1,...,n} \| T(j+1,i;x)-T(j,i;x)\| \le $$
$$ \max \{ \| A(j;x)- A(j-1;x) \| t', \mbox{ }
\| A(j;x)^{j+1}A(j+1;x)/(j+1)!\| \} \le \| P\| p^z,$$
where $z=-(j+1)$, $t'=\| P\| $, since the second term in $\{ , \} $ is
less than the first and 
$$\| T(j+1,i;x)-T(j,i;x)\| =\| P(i;x)-P(i;x)+ 
\{ \sum_{s=1}^{j-1} A(j-1;x)^s(T(j+1,i;x)- $$
$$T(j,i;x))/s! \} + \{ \sum_{s=1}^{j-1} 
(A(j;x)^s -A(j-1;x)^s)T(j+1,i;x)/s! \} $$
$$+ A(j;x)^{j+1}T(j+1,i,;x)/(j+1)! \| .$$
Therefore, there exists a sequence satisfying Formulas
$(i-iii)$ for each $j$. Hence there exists 
$$\lim_{j\to \infty } A(j;x)=A(x)$$ such that
$A: C_0(\tau , M\to X) \to C_0(\tau ',M\to X)$. This mapping may be
considered as a vector field on $M$ of class $C_0(\tau )$, $A(x)\in Vect_0
(\tau ,M)$, consequently, there exists 
$$\lim_{j\to \infty } g(j;x)=g(x)\in C_0(t,M\to X).$$
In view of $\| A(j;x)^{s+j} /(s+j)!\| _{\tau }
\le p^z$, where z$=-2(j+s)+(j+s)/(p-1)$ for each $s \in {\bf N}$
there is 
$$exp \{ qA(j;x) \} x=g^q(j;x)\mbox{ with }\lim_{j\to \infty }
g^q(j;x)=g^q(x)\in W$$ 
(that is convergent relative to $\rho _0^{\tau }$)
for each $g(x) \in W$ and $q \in B({\bf K},0,1)$ such that $g^1(x)=g(x)$.
Moreover, to  $\{ g^q(x): q \in B({\bf K},0,1) \}$
there corresponds a one-parameter subgroup in $W$, where $q \in {\bf Z_p}$,
since $\| qA(j;x)\|_{\tau } \le p^{-2}$ for each $q, y\in B({\bf K},0,1)$.
\par  Indeed, $g^q=g_{i,j}^q$ are given as mappings from $\phi _j(U_j)$
into $\phi _i(U_i)$ for a given $i, j$,
$\| g_{i,j}^q-id_{i,j} \| _{\tau } \le 1/p$, so $g_{i,j}^q$
generate $g^q \in W$, $g^q: M\to M$, since $g^q$ is an isometry,
consequently, $g^q: \in G(\tau ,M)$. For $t=
\infty $ we consider all $\tau \in \bf N$.
\par  In general, for each $f \in G(t,M)$ there is a sequence
$\{ f_l(x): l {\in \bf N} \} \subset
G(t,M)$ such that in local coordinates $x=\{ x(i) : i \in I \}
\in B(c_0(I,K),0,r)$ for each $i >l$ 
the following condition is satisfied $(f^i_l(x))=x(i)$
and there exists $A_l(x) \in Vect_0(\tau ,M)$ with the corresponding
$g_l^q(x)\in W$ and $g^1_l(x)=f_l(x)$ for each $x \in M$,
where $\lim _{j \to \infty } \| f_l -f\| _{\tau }=0$.
Then
$$\lim _{l \to \infty} g_l^q(x)=\lim_{l\to \infty }g^q(x)\in W$$ 
converges
relative to $\rho _0^{\tau }$ and $A(x)=\lim_{l\to \infty } A_l(x)$
with $\| A\| _{\tau }\le p^{-2}$, where 
$$A=\sum_{m, i} a(m,A^i){\bar Q}_m(x) \partial _i \in Vect_0(\tau ,M),$$ 
$a(*,*) \in K$, that is, for each $c>0$
the set $\{ (i, Ord(m), |m|):$ $|a(m,A^i)| J(\tau ,m)>c \}$ is finite.
\par  The field $\bf K$ is equal to the disjoint union $\bigcup_{j\in \bf N}
B({\bf K},k_j,1)$, where $k_j \in \bf K$, $k_1=0$. Defining $g^{q
+k_j}=g^q$ for $j>1$ and $q \in B({\bf K},0,1)$, we get the extension
of class $C_0(\tau )$ by $q$
for $g^q$ from $B({\bf K},0,1)$ onto $\bf K$ by $q$,
for $1\le t\le \infty $. For $t=an_r$ we use
the additive group $B({\bf K},0,1)$.
Then $\partial g^q(x)/\partial q=A(x)g^q(x)$
for each $q \in B({\bf K},0,1)$ and $x \in M$,
$A=\sum_iA^i\partial _i$, $A^i\in C_0(\tau ,M\to {\bf K})$.
\par In the cases of the non-local filed $\bf F$ or
$G=Diff(t,M)$ consider the family $\Upsilon =\{ 
\eta _{z_1,...,z_n,\bf K} \} $
of all embeddings $\eta _{z_1,...,z_n,\bf K}: 
sp_{\bf K} \{ z_1,...,z_n \} \hookrightarrow X$, where 
$\bf K$ are all possible local subfields of $\bf F$
and $z_1,...,z_n$ are linearly independent vectors in $X$,
$n\in \bf N$. If $f\in G$, then $f: M_{z_1,...,z_n,\bf K}
\to f(M_{z_1,...,z_n,\bf K})$ is the diffeomorphism
of class $C_0(t)$, where $M_{z_1,...,z_n,\bf K}:=
M\cap \eta _{z_1,...,z_n,\bf K} (sp_{\bf K} \{ z_1,...,z_n \} ).$
Let $\rho ^{\tau }$ be a left-invariant ultrametric in $G$ induced by 
the norm in $C(\tau ,M\to X)$ for $M$ with the finite atlas.
There are embeddings of spaces $G(t,M_{z_1,...,z_n,\bf K})$
into $G$ such that $\rho ^{\tau }$ induces the equivalent 
ultrametric $\rho ^{\tau }_0$ in $G(t,M_{z_1,...,z_n,\bf K})$.
Therefore, there exists a clopen subgroup $W$ in $G$ such that
for each $f\in W$ and its restriction $f|_{M_{z_1,...,z_n,\bf K}}$
there exists a one-parameter family $\{ g^q_{z_1,...,z_n,\bf K}
: q\in {\bf K} \} $ which has an embedding into 
$W|_{M_{z_1,...,z_n,\bf K}}$. These families can be chosen 
consistent on $M_{z_1,...,z_n,\bf K}\cap M_{y_1,...,y_m,\bf L}$,
since ${\bf K}\cap \bf L$ is a local field for two local fields
$\bf K$ and $\bf L$ such that ${\bf Q_p} \subset {\bf K}\cap {\bf L}
\subset {\bf K}\cup {\bf L}\subset \bf F$, moreover, there exists
a local field $\bf J$ such that ${\bf K}\cup {\bf L}
\subset \bf J$. 
This means, that $g^q_{z_1,...,z_n,\bf K}(x)
=g^q_{y_1,...,y_m,\bf L}(x)$ for each $x\in M_{z_1,...,z_n,\bf K}
\cap M_{y_1,...,y_m,\bf L}$ and for each $q\in {\bf K}\cap \bf L$.
Hence there exists $g^q(x)$ for each $x\in 
cl_M \{ \bigcup_{z_1,...,z_n, \bf K} M_{z_1,...,z_n,\bf K} \} $
and each $q\in cl_{\bf F}\{ \bigcup_{{\bf K}\subset \bf F}
{\bf F} \} $, where $cl_MA$ denotes a closure of a subset $A$ in $M$.
In $\bf C_p$ the union of all local subfields is dense.
If $\bf F$ is not contained in $\bf C_p$, then it can be constructed
from $\bf C_p$ with the help of operations of spherical completion
$\bf C_p^U$ or of quotients of definite algebras over
$\bf C_p$ or $\bf C_p^U$ and so on by induction \cite{diar}.
Therefore, these consisitent 
families generate a one-parameter subgroup $\{ g^q: q\in B({\bf F},0,1) \} $
in $W$ such that $g^1=f$.
\par  Let now $M$ be with a countable infinite atlas and 
at first $1\le t\le \infty $ then
from the definition of topology in $G$ the following set
$$W:=\{ f\in G: \mbox{ } supp(f)\subset U^{E(f)},
\mbox{ } E(f)\subset {\bf N}, card (E(f))<\aleph _0,
\rho _{0,U^{E(f)}} ^{\tau }(f,id)\le p^{-2} \} $$ 
is the clopen subgroup, 
where $\rho _{0,U^E}^{\tau }(f,g)$ are ultrametrics in $G(t,U^E)$
inducing pseudoultrametrics in $G$,
$\tau =t$ for $t\ne \infty $ and $1\le \tau <\infty $ for $t=\infty $,
$U^E:=\bigcup_{i\in E}U_i$, $(U_i,\phi _i)$ are charts of $M$.
\par  For $t=la $ let 
$$W:= \{ f\in G: \mbox{ } supp(f)\subset U^{E(f)},
E(f)\subset {\bf N}, card (E(f))<\aleph _0,$$
$$\rho _{0,U^{E(f)}}^{an_r}(f,id)\le p^{-2}, f\in C_0(an_r,
M\to M), r\in \Gamma _{\bf F} \} , $$
where $\rho _{0,U^E}^{an_r}(f,g):=\sup_{i\in E}\|(g^{-1}\circ f-
id)_{i,j}\| _{an,r,E}$. 
\par  $\bf (2).$ Let at first $t=an_1$. Let us prove that
the function $exp :Vect(t,M)\to Diff(t,M)$ is not locally
bijective. Let $M=B({\bf F} ,0,1)$ be a manifold over $\bf F$. 
We suppose, that there exists $q\in \bf F$ such that $q^l\ne 1$ for
each $l=1,...,n-1$, $q^n=1$, where $n$ is not divisible by $p$
and $1<n\in \bf N$,
$q^s \in \bf F$ and $\mid q^s\mid  _p=1$ for each $s \in {\bf Z_p}
\cap \bf F$. Further $q^sM=M$ (acts as
the multiplication $x\mapsto q^sx$ for each $x \in M$) 
and $q^s \in Diff(t,M)$, particularly,
for $s=1$, $(q^1)^n =id$. Let $H:= \{ g: g \in Diff(t,M), g^n =id \}$,
consequently, $gq^1g^{-1} =q^1=q$ for each $g \in H$. Whence $q^1$
belongs to each one-parameter subgroup $gTg^{-1}$ in $Diff(t,M)$, where
$T:= \{ q^s: s \in B({\bf F},0,1) \}$.
\par Now we consider the case, when the field $\bf F$ has not 
sufficient roots of unity. If $G$ would be a Banach-Lie group,
then there will exists a clopen subgroup $W$ in $G$ such that
the Campbell-Hausdorff formula \cite{boug} will be valid. 
Let $g^q_m=exp(qx^m\partial )x=$ $\sum_{l=0}^{\infty }(qx^m\partial )^nx/n!$,
where $x\in M=B({\bf F},0,1),$ $0\le m\in \bf Z$, $q\in B({\bf F},0,1/p)$. 
Therefore, $g^q_0(x)=x+q$, $g^q_m(x)=\sum_{k=0}^{\infty }
q^kx^{k(m-1)+1}\gamma (k,m)/k!$, where
$\gamma (k,m):=m(2m-1)(3m-2)...((k-1)m-k+2)$ for $k\ge 1$,
$\gamma (0,m)=1$, $0!=1$. Then 
$$(ad \mbox{ }u)^s(v)=
\xi ^s\zeta x^{n+s(m-1)}(n-m)(n-1)(n+m-2)...(n+(s-2)m-s+1)$$
for each elements $u=\xi x^m\partial $ and 
$v=\zeta x^n\partial \in {\sf g}:=T_eG$, where $\xi , \zeta
\in \bf F$, $u:=\xi \partial $ for $m=0$.
Let $w=ln(e^u\circ e^v)$ be given by the Campbell-Hausdorff formula.
Then calculating several lower terms of the series
we get that $e^w(x)$ does not coincide with $g^{\xi }_n\circ 
g^{\zeta }_n(x),$ where $\xi , \zeta \in B({\bf F},0,1/p)$.
This contradicts our supposition, consequently, $G$ is not the 
Banach-Lie group. For $dim_{\bf F}X>1$ it is sufficient to consider
embeddings of $Diff(an_1,B({\bf F},0,1))^k$ into $Diff(an_r,M)$,
where $1<k\le dim_{\bf F}X$. 
\par In the case of $0\le t\le \infty $ for each $f\in W$
there exists an infinite family $g^q_l$ of one-parameter subgroups
such that $g^1_l=f$ and $\partial g^q_l/\partial q
=\partial g^q_i/\partial q$ for each $i, l \in \bf N$, 
$q\in B({\bf F},0,1),$ since we can consider locally-constant 
additional terms for a given $g^q$.
\par Each subgroup $G(t,U^E)$ for $1\le t\le \infty $ or
$G(an_r,U_j)$ for $\phi _j(U_j)\subset B(X,\phi _j(x),r)$
are closed in $G$ and are not the Banach-Lie groups ,
consequently, $G$ is not the Banach-Lie group.
\par {\bf 3.2. Theorem.} {\it Let groups $G:=Diff(t,M)$ and
$G:=G(t,M)$ be given by Definition 2.4.
Then $G$ is simple and perfect.}
\par  {\bf Proof.} It is sufficient to consider the case
of a manifold $M$ on a complete locally $\bf F$-convex
space $X$, since the perfectness and simplicity of
$G$ and its completion $\tilde G$ are equivalent. 
Consider at first $G$ with
$t\ge 1$ or $t=an_r$. If $f, g \in W\subset G(t,M)$ (see Theorem 3.1), then 
there are vector fields $A_f$ and $A_g$
of class $C_0(t)$ on $M$ and one-parameter subgroups
$f^q$, $g^q\subset W$, $q \in B({\bf F},0,1)$ such that $\partial
f^q/\partial q=A_ff^q$ and $\partial g^q/\partial q=A_gg^q$, where
$A_f(x)=A^i_f(x)\partial _i$, the summation is accomplished by $i \in I$,
$I$ is an ordinal. Let $A^i$ be of class $C_0(\tau )$ with
$\tau =\infty $ for $1\le t\le
\infty $ or $\tau =t$ for $t=an_r$, then elements $exp(qA(x))x$
are dense in $W$ for such $t$, where $A=A^i\partial _i$,
$q\in B({\bf K},0,1)$. For $B=\bar A^i\partial _i$ with $\bar A^i=
\delta _{i,j}$, where $j\in I$ is fixed, $\delta _{i,j}=1$ for
$i=j$ and $\delta _{i,j}=0$ for $i\ne j \in I$, $W \ni exp(qB)x\ne id$.
If $C\in Vect_0(\tau ,M)$ with $\tau =\infty $ or $\tau =an_r$,
also $\| C\|_{\tau '}\le p^{-2}$ ($\tau '\in \bf N$ or $\tau '=an_r$
respectively), then there exists $A\in Vect_0(\tau ,M)$ with
$\| A\|_{\tau '}\le p^{-2}$ such that their commutator $[A,B]=C$.
Indeed, $[A,B]=A^i(\partial _iB^k\delta _{k,j})\partial _j$ $-B^k\delta
_{k,j}(\partial _jA^i)\partial _i$ $=-(\partial _jA^i)\partial _i=
C^i\partial _i$. In view of \S 79.1 and \S 80.3 \cite{sch1} there is
the antidifferentiation
$P(x^j)$ by the variable $x^j$ (in the local coordinates of $At(M)$)
such that $A^i(x)=-(P(x^j)C^i)(x)$. From this it follows that $[W,W]=W$,
where $[W,W]$ is the minimal subgroup in $W$ generated by the following
subset $\{ [g,f]:=g\circ f\circ
g^{-1}\circ f^{-1}| g, f \in W\} $. Therefore, $W$ is perfect. 
\par  Suppose there is a normal subgroup $V$ in $W$, $V\ne \{ e\} $ and
$V\ne W$, then $gVg^{-1}=V$ for each $g \in W$.
Let $v\subset w$ be corresponding to $V$ and $W$ subsets in the algebra
$Vect_0(\tau ,M)$, hence $[v,w]\subset v$, where $[A,B]$ is the
commutator in the algebra $Vect_0(\tau ,M)$. Therefore, there are
$A\in w\setminus v$
and $0\ne B\in v$. Since $[v,w]\subset v$, then $[p^2\partial _i,
B]\in v$ for each $i$, so it may be assumed that there is $i\in I$
with $a(0,B^i)\ne 0$.
\par  For $Vect_0(an_r,M)$ we get the equations
$\sum_{i, m+n=k+e_i}(b(n,C^i)b(m-e_i,B^j)-$ $b(m,B^i)b(n-e_i,C^j))=$
$b(k,A^j)$, consequently, solving them recurrently we find $B\in v$
and $C\in w$ for which $[C,B]
=A$. This is possible, since for each $c_l=p^{-l}$, $l\in \bf N$,
the family $\{ (i,|n|,Ord(n)): |b(n,A^i)|r^{|n|}>c_l \} $ is finite
for $A \in Vect_0(an_r,M)$, where $b(m,B^j)\in \bf F$ are 
expansion coefficients by $x^m:=x_1^{m_1}...x_n^{m_n}$, $x=(x_1,...,x_n)$,
$x_i\in \bf F$, $m=(m_1,...,m_n)$, $0\le m_i\in \bf Z$, $n\in \bf N$.
\par  Locally analytic functions \cite{sch1} are dense in
$C_0(t,M\to X)$ for $1\le t\le \infty $, hence, $[v,w]=w$
contradicting our assumption. Therefore, $W$ is simple; that is,
it does not contain any normal subgroup $V$ with $V\ne \{ e\} $ and
$V\ne W$ at the same time.
\par  The group $G$ is the disjoint union of $g_iW$, $G=
\bigcup_{j\in \bf N}g_jW$, such that $\rho _0^t(g_j,
g_k)>p^{-2}$ for $j\ne k$, hence for chosen $\{ g_j:
j\in {\bf N} \} $ with $g_1=e=id$ and each $f\in G\setminus W$
there is the unique $j$ and $f_2\in W$ with $f=g_jf_2$. From $\bar
Q_m(x+c)=\bar Q_m(x)$ for $|m|>0$ and each $x \in B({\bf K},0,1)$
it follows that $c=0$, where $n=Ord(m),$ $\bar Q_m$ are basic 
Amice polynomials \cite{ami,lusmb98}. 
Then considering all $g\in W$
having the form $(id+\xi \bar Q_m(x)e_i)$ in local coordinates with
$\xi \in B({\bf K},0,p^{-2})$
we get $[G\setminus W,W]\supset G\setminus W$ and $\{ gfg^{-1}:
g \in W \} $ $\ne \{ f\} $ for each $f \in G\setminus W$,
hence $[G,G]=G$.
\par  Suppose there is a normal subgroup $V\subset G$, $V\ne e$,
$V\ne G$. Then for each $f$ and $g\in V$ with $fg^{-1}\in W$
we get that $fg^{-1}=e$, since $V\cap W$ is the normal subgroup
in $W$, consequently, $V$ is discrete and $\rho _0^t(f,g)>p^{-2}$
for each $f\ne g \in V$. Therefore, $hfh^{-1}=f$ for each $h \in W$,
but this contradicts the statements given above. Therefore,
$G$ is simple and perfect.
\par  For $0\le t<1$ the group $G(1,M)$ is dense in $G(t,M)$,
for $t=la$ we use the inductive limit topology,
consequently, $G(t,M)$ also is simple and perfect.
The case of $Diff(t,M)$ follows from the case of $G(t,M)$ 
analogously to the proof of Theorem 3.1.
\par  {\bf 3.3. Notes.} For each
chart $(U_i, \phi _i)$ there exists
a tangent vector bundle $TU_i=U_i\times X$, 
consequently, $TM$ has the following
atlas $At(TM)=$ $\{ (U_i\times X,\phi _i\times I): i \in \Lambda \} $,
where $I: X\to X$ is the unit mapping, $\Lambda \subset \bf N$, $TM$
is the tanget vector bundle over $M$.
\par  Suppose $V$ is a vector field on $M$. 
Then by analogy with the classical
case we can define the following mapping $\bar exp_x(zV)=x+zV(x)$, 
where for the corresponding pseudoultranorm $u$ in $X$
and sufficiently small $\epsilon >0$ from $u(V(x))|z|\le \epsilon $ 
it follows $x+zV(x)\in U_i$ for each $x \in U_i$ such that
$\phi _i(x)$ is also denoted by $x$, $z\in \bf K$. 
Moreover, there exists a refinement
$At'(M)=$ $\{ (U'_i,\phi '_i): i\in \Omega \} $ of $At(M)$
such that $\phi '_i(U_i)$ are $\bf F$-convex in $X$. 
The latter means that $\lambda x+(1 - \lambda )y\in \phi '_i(U'_i)$
for each $x, y \in \phi '_i(U'_i)$ and each $\lambda \in B({\bf F},0,1)$.
The atlas $At'(M)$
can be chosen such that $(\bar exp_x|_{U'_i}): x\times \{
\phi _i(U'_i) - \phi _i(x) \} 
\to U'_i$ to be the analytic isomorphism for each
$i \in \Omega $, where $x\in U_i'$.
\par  Then $T_fC_*(t,M\to N)=$ $\{ g
\in C(t,M\to TN)|$ $\pi \circ g=f \} $, consequently, $C_*(t,M\to TN)=$
$\bigcup_{f \in C_*(t,M\to N)}T_fC_*(t,M\to N)$, 
where $\pi : TN \to N$ is the natural projection.
Therefore, the following mapping $\omega _{\bar exp}:$
$T_fC_*(t,M\to N) \to $ $C_*(t,M\to N)$ such that 
$\omega _{\bar exp}(g)={\bar exp}\circ g$ is defined. This gives charts
on $C_*(t,M\to N)$ induced by charts in $C_*(t,M\to TN)$.
\par {\bf 3.4. Theorems.} {\it Let $G=Diff(t,M)$ and
$G=G(t,M)$ be the same as in \S 2.4, where
$1\le t \le \infty $ or $t=la$, $\bf F$ and $X$ are complete. Then 
\par $\bf (1)$ if $V$ is a $C(t)$-vector field on $M$,
then its flow $\eta _z$
is a one-parameter subgroup by $z \in B({\bf F},0,1)$ in $G$; 
\par $\bf (2)$ the mapping $z\mapsto \eta _z$ is continuously differentiable
by $z\in B({\bf F},0,1)$; 
\par $\bf (3)$ the mapping $\tilde Exp: T_{id}G\mapsto G$,
$V\mapsto \eta _z$, is continuous and defined in a neighbourhood
of $0$ in $T_{id}G$ for each $z\in B({\bf F},0,1)$, 
the mapping $(z,V)\mapsto
\eta ^V_z$ is of class $C(t)$;
\par $\bf (4)$ $G$ is an analytic manifold and for it
the mapping $\tilde E:$
$TG\to G$ is defined such that $\tilde E_{\eta }(V)=\bar exp_{\eta (x)}
\circ V_{\eta }$ from some neighbourhood  $\bar V_{\eta }$
of the zero section in $T_{\eta }G\subset TG$
onto some neighbourhood $W_{\eta }\ni id \in G$, 
such that $\bar V_{\eta }
=\bar V_{id}\circ \eta $, $W_{\eta }=W_{id}\circ \eta $,
$\eta \in G$ and $\tilde E$
belongs to the class $C(\infty )$ by $V$, 
$\tilde E$ is the uniform isomorphism
of uniform spaces $\bar V$ and $W$.}
\par  {\bf Proof.} As in the proof of Theorem 3.1
we reduce the consideration to the case of $M$ with a finite atlas
on the Banach space $X$ over $\bf F$
and $G=G(t,M)$ and then generalize results for infinite atlases 
on the locally $\bf F$-convex space $X$
and $G=Diff(t,M)$ using inductive limits of spaces $C(t,U^E\to Y)$
and the projective limit $X=pr-\lim X_u$.
\par For each submanifold $M_n$
in $M$ with the embedding $\theta : M_n\hookrightarrow M$ and
$dim_{\bf F}M_n=n$
let us consider $V|_{M_n}: M_n\to TM$, $\pi \circ V(x)=x$.
Therefore, in view of Theorem 6.1 \cite{khr} 
(and analogously we get existence of solutions in classes $C(t)$)
there is the
solution $\eta _z$ for some $c>0$ and all $z \in B({\bf F},0,c)$,
that is, $\partial \eta _z(x)/\partial
z=V(x)\eta _z(x)$, $\eta _0(x)=x$ are dependent upon $x \in M$,
$\eta _0=id$,
$\eta ^V_z(x)=\eta _z(x)$ are dependent upon $V$.
This local solution is unique in the analytic case, but it is not unique
in $C(la)$ and $C(t)$ classes.
Here a constant $\infty >c>0$ depends only upon $0<R<\infty $, 
$M$ and $t$, where $V$ is in the neighbourhood of the zero section
$B(T_{id}C_*(\tau ,
M\to M),0,R)$ and the ultranorm on the Banach space $T_{id}
C_*(t,M\to M)$ is inherited from the Banach space $C_*(t,M\to TM)$.
\par  The clopen ball $B({\bf K},0,c)$ is the additive subgroup in
$\bf K$, consequently, $z\mapsto \eta _z$ is the homomorphism
such that $z_1+z_2\mapsto
\eta _{z_1+z_2}=\eta _{z_1}\circ \eta _{z_2}$, $\eta _0=id $. Moreover,
$z\mapsto \eta _z$ and $V\mapsto
\eta _z=\eta ^V_z$ are $C^{\infty }$-mappings by $V$ 
and $z$ in some neighbourhoods of $0$. 
On the other hand, $B({\bf F},0,1)$ is a disjoint union 
of balls of radius $0<c<1$. Hence there exists a solution for each
$z\in B({\bf F},0,1)$ (see aslo \S 45 in \cite{sch1}).
\par Then there are $R$ and $c$ such that 
$\rho _*^t(\eta ^V_z,id)
\le 1/p$ for each $z \in B({\bf F},0,c)$ and $V \in B(T_{id}C_*(t,M\to M),
0,R)$, hence for $Rc=R'c'$, $c'=1$, we get the following mapping
$V\mapsto \eta ^V_z$ for each 
$V \in B(T_{id}C_*(t,M\to M),0,R')$, where $z\in
B({\bf F},0,1)$. Then $V\mapsto \eta _1$ generates the mapping
$\tilde Exp(V)=\eta _1$.
Hence, $\tilde Exp$ is defined in the neighbourhood of $0$ in
$T_{id}G$ and $\tilde Exp \in C^{\theta }(\infty ,B(T_{id}G,0,R')\to G)$,
where the last space is given relative to the mapping $\theta
=\tilde \pi _{id}: T_{id}G\to G$ being the natural projection.
\par Let $V(\eta )\in T_{\eta }G$ for each $\eta \in G$ and
$V\in C_*(t,G\to TG)$, where $\tilde \pi \circ V(\eta )=\eta $,
$\tilde \pi : TG \to G$, $V$ is a vector field on $G$ of class
$C_*(t)$. If $\tilde V:= \{
g \in C_*(t,M\to M):$ $\rho _*^{\kappa }(g,id)\le 1/p \} $
and $g \in \tilde V$, where $\kappa =t$ for $t \ne \infty $ 
and ${\bf N}\ni \kappa \ge 1$
for $t=\infty $, then $g: M\to M$ is
the isometry, consequently,
$g\in Hom(M)$, that is, $\tilde V\subset G$ and $G$ is a neighbourhood of
$id$ in $C_*(t,M
\to M)$. Since $M=\bigcup_{i\in \Lambda }U_i$, $TM=\bigcup_i
(U_i\times X)$, then $C_*(t,M\to M)$ and $C_*(t, M\to TM)$
have atlases with clopen charts. The $C^{\infty }$-atlases 
$At(C_*(t,M\to M))$ and  $At(C_*(t,
M\to TM)$ induce clopen atlases in $G$ and
$TG$, since $M$ and $\bar exp$ are 
of class not less than $C^{\infty }$ (see \S 2.4 and \S 3.3). 
\par  The right multiplication $R_f: G\to G$, $g\mapsto g\circ f=R_f(g)$ and
$\alpha _h: C_*(t,M'\to N)\to C_*(t,M\to N)$, $\alpha _h
(\zeta )=\zeta \circ h$ for $f \in G$ and $h \in C_*(t,M\to M')$
belong to the class $C(\infty )$ for $1\le t \le \infty $, 
also to $C(an_r)$ for $t=an_r$ (see Theorems 2.6). Let $g\in
\tilde V$, then $g=id+Y$ with $\| \phi _i(Y|_{U_j}) \| _{
C_*(t,U_j\to X)}
\le 1/p$ for each $j$ (see \S 2.4), hence, $\tilde g_z=
id+zY\in \tilde V$ for each $z\in B({\bf F},0,1)$ and $(\partial R_f
\tilde g_z/\partial z)|_{z=0}=R_fY$. From this it
follows that
each vector field $V$ of class $C_*(t)$ on $G$ is right-invariant,
that is, $R_fV_{\eta }=V_{\eta \circ f}$ for each $f$ and $\eta \in G$,
since $G$ acts on the right on $G$ freely and transitively (that is,
$g\circ f=f$ is equivalent to $f=id$, $Gf=fG=G$). 
\par Therefore, the vector field $V$ on $G$ of class $C_*(t)$ has the
form $V_{\eta (x)}=v\circ \eta (x)=v(\eta (x))$, where $v$ is a vector
field on $M$ of the class
$C_*(t)$, $\eta \in G$, $v(x)=V(id(x))$. Since
$\bar exp: TM\to M$ is locally analytic on corresponding charts, then
$\tilde E(V)=\bar exp\circ V$ has the necessary properties (see
for comparison also the classical case in 
\cite{bao} and in \S 3, \S 9 \cite{ebi}).  
\par {\bf 3.5. Notation.} Let the group $G=G(t,M)$ be given by
\S 2.4, where ${\bf Q_p}\subset {\bf F} \subset \bf C_p$,
an atlas of $M$ is countable.
A complete locally $\bf F$-convex space $X$ has a structure of a locally 
$\bf K$-convex space $X_{\bf K}$ over a local subfield $\bf K$ 
in $\bf F$, then $X_{u,\bf K}=X_{\bf K}/ker (u)$ is isomorphic with
the Banach space $X_{u,\bf K}=c_0(J_u,{\bf K})$ 
over a local field $\bf K$, where $J_u$ is an ordinal, 
$u$ is a pseudoultranorm in $X$ (see \cite{nari}, Ch. 5 in \cite{roo},
\cite{lusmb98}). 
There exists $M_{\bf K}$ which is a manifold $M$ 
on $X_{\bf K}$ instead of $X$. The projection 
$\pi _u: X_{\bf K}\to X_{u,\bf K}$ induces a projection
$\pi _u: M_{\bf K}\to M_{u,\bf K}$, where $M_{u,\bf K}$
is a manifold on $X_{u,\bf K}$.
Let each $X_u:=X/ker (u)$ be of separable type over $\bf F$
for a family of pseudoultranorms $ \{ u \} $ defining topology of $X$.
Let us consider the following space
$$C_{0,a}(t,M_{u,\bf K}\to X_{u,\bf K}):= \{ f\in C_0(t,M_{u,\bf K}
\to X_{u,\bf K}): \mbox{ } \| f\|_{t,a}:=$$
$$\sup_{k,i,j,m}
[|a(m,\phi _k\circ f^i|_{U_j})| J_j(t,m)\max (p^i,p^{Ord(m)+|m|})]<
\infty ,$$ 
$$\lim_{j+i+|m|+Ord(m)\to \infty }|a(m,f^i_{U_j})|J_j(t,m)\max (
p^i, p^{|m|+Ord(m)})=0 \} $$ 
for $t\ne \infty $, $C_{0,a}(\infty ,M_{u,\bf K}\to X_{u,\bf K}):=
\bigcap_{l\in \bf N}C_{0,a}(l,M_{u,\bf K}\to X_{u,\bf K})$,
where $(U_j,\phi _j)$ are charts of $At(M_{u,\bf K})$ 
with omitted index $(u,{\bf K}),$
$J_j(t,m):=$ \\ 
$ \| {\bar Q}_m|_S \| _{C_0(t,S\to X_{u,\bf K})},$
$S:=(U_j)_{u,\bf K}\cap sp_{\bf K} \{ e_1,...,e_{Ord(m)} \} $,
$\{ e_j: j \} $ is the standard orthonormal base of
$c_0(J_u,{\bf K})$.
\par For the manifolds $M$ and $N$
with a given mapping $\theta : M \to N$ (see \S 2.4) 
we define an ultrauniform space
$$C_{0,a}^{\theta }(t,M_{u,\bf K}\to N_{v,\bf K}):=\{ f \in 
C_0^{\theta }(t,M_{u,\bf K}\to N_{v,\bf K})|
(f_{i,j}-\theta _{i,j}) \in C_{0,a}(t,\phi _j(U_j)\to Y
_{v,\bf K})$$ 
$$\mbox{for each }i \mbox{ and }j, \mbox{ where } 
\rho ^t_a(f,\theta ):=\sum_{i,j}\|
(f-\theta )_{i,j}\| _{C_{0,a}(t,\phi _j(U_j)\to Y_{v,\bf K})}<\infty \} $$
for each $0\le t <\infty $. There exists a subgroup
$$G_a(t,M_{u,\bf K}):=
G(t,M_{u,\bf K})\cap C_{0,a}^{id}(t,M_{u,\bf K}\to M_{u,\bf K}) $$
with an ultrametric $\rho _a^t(f,id)$
for $\theta =id$ and $0\le t< \infty $.
\par {\bf 3.6. Theorems.} {\it Let 
$X,$ $\bf F$, $G:=G(t,M)$ and $G_a(t,M_{u,\bf K})$ be given
by \S 2.4 and \S 3.5. Then
\par $\bf (1).$ For each $0\le t\le \infty $
spaces $G_a(t,M_{u,\bf K})$ and
$C_{0,a}(t,M_{u,\bf K}\to X_{u,\bf K})$ are separable and complete. 
\par $\bf (2).$ Each $G_a(t,M_{u,\bf K})$ is $\sigma $-compact and 
$G_a^r(t,M_{u,\bf K}):=B(G_a(t,M_{u,\bf K}),id,r):= 
\{ g\in G_a(t,M_{u,\bf K})|$
$\rho _a^t(g,id)\le r \} $ has an embedding as 
a compact separable subgroup in
$G(t,M)$ in the topology inherited from it for $0\le t< \infty $
and $0<r<\infty $.
\par $\bf (3)$ $T_eG(t,M)\subset Vect_0(t,M)$ (see \S 3.2), 
moreover, $sp_{\bf K}\bigcup_u T_eG_a(t,M_{u,\bf K})$ and
$sp_{\bf K}\bigcup_u \mbox{ } T_eB(G_a(t,M_{u,\bf K}),id,r)$ are contained in
$T_eG(t,M_{\bf K})$ and dense in it for $1\le t\le \infty $.
\par $\bf (4).$ In $G$ for $0\le t<\infty $ there is a family
$\{ G^n_{u,\bf K}: n\in {\bf N}, u, 
{\bf Q_p}\subset {\bf K}\subset {\bf F} \} $ 
of compact subgroups such that $G^n_{u,\bf K}
\subset G^{n+1}_{u,\bf K}$ for each $n$ 
and each local subfield $\bf K$ in $\bf F$, moreover,
$\bigcup_{n\in {\bf N}, u, {\bf K} }G^n_{u,\bf K}=:
\bar G_a(t,M)$ is dense in $G$.}
\par  {\bf Proof.} From Formulas $2.6(i,ii)$ it follows that
$G_a(t,M_{u,\bf K})$ are the complete topological groups and $C_{0,a}
(t,M_{u,\bf K}\to X_{u,\bf K})$ is the complete locally $\bf K$-convex space
(and it is the Banach space for $0\le t< \infty $). 
They are separable and Lindel\"of, since
$M_{u,\bf K}$ and $X_{u,\bf K}$ are separable. 
\par The uniformity in $G_a(t,M_{u,\bf K})$ is given by 
the right-invariant ultarmetric
$\rho _a^t(f,g)=\rho _a^t(g^{-1}f,id)$ for $t\ne \infty $
and by their family
$\{ \rho _a^l: l \in {\bf N} \} $ for $t=\infty $, where
$f$ and $g\in G_a(t,M_{u,\bf K}).$
Therefore, $B(G_a(t,M_{u,\bf K}),id,r)=:G^r_{u,\bf K}$ 
is also the topological group, which
is clopen in $G_a(t,M_{u,\bf K})$. 
The Banach space $C_0(t,M_{u,\bf K}\to X_{u,\bf K})$ 
is linearly topologically
isomorphic with $c_0(\bf \omega _0,{\bf K})$ and subsets 
$S:= \{ x\in c_0(
\omega _0,{\bf K}):$ $|x(i)|_{\bf K}\le p^{-k(i)}\mbox{ for each }
i\in {\bf N} \} $ are sequentially compact in $c_0(\omega _0, {\bf K})$ for
$\lim_{i\to \infty }k(i)=\infty $, consequently, $S$ are compact
\cite{eng}. In addition
$sp_{\bf K}S$ is dense in $c_0(\omega _0,{\bf K})$. 
\par Analogously to the proof of Theorems 3.4 we get neighbourhoods
$\tilde U\ni 0$ in $T_{id}G_a(t,M_{u,\bf K})$ and $\tilde W\ni id$ in
$G_a(t,M_{u,\bf K})$, such that $\tilde E|_{\tilde U}:$ 
$\tilde U\to \tilde W$ is the uniform isomorphism. 
There exists an embedding of $G^r_{u,\bf K}$ into $G(t,M)$, 
since each function 
$f$ on $M_{u,\bf K}$ has an extension $\tilde f$ on $M_{\bf K}$ such that 
${\tilde f}|_{M_{\bf K}\ominus M_{u,\bf K}}=id$, where
the decomposition $M_{\bf K}=
(M_{\bf K}\ominus M_{u,\bf K})\oplus M_{u,\bf K}$ is induced by the 
projection $\pi _u$, since $\bf K$ is spherically complete.
Due to $\tilde W\supset G^r_{u,\bf K}$ for
$0<r\le p^{-2}$ the subgroup $G^r_{u,\bf K}$ is compact in 
the weaker topology inherited from $G(t,M)$. For
$p^{-2}< r<\infty $ considering in local cooordinates basic functions
${\bar Q}_me_i$ and using the ultrametric $\rho _a^t$ we get for
$f \in G^r_{u,\bf K}$ that only a finite number of $(m,k,i,j)$ are such that
$|a(m,\phi _k\circ f^i|_{U_j})|>p^{-2}$.
Therefore, $G^r_{u,\bf K}$ is compact. 
From $G_a(t,M_{u,\bf K})=\bigcup_{i\in \bf N}
g_iG^r_{u,\bf K}$ for some family $ \{ g_i:$ $i \in {\bf N} \} 
\subset G_a(t,M_{u,\bf K}) $
(or $G_a(t,M_{u,\bf K})=\bigcup_{l\in \bf N}G^l_{u,\bf K}$) 
it follows, that $G_a(t,M_{u,\bf K})$ is $\sigma $-compact
in $G(t,M)$. 
\par  In view of \S 3.2 and \S 3.4 we get
$T_eG(t,M)\subset Vect_0(t,M)$ and $T_eG^r_{u,\bf K}
\subset T_eG_a(t,M_{u,\bf K})\subset T_eG(t,M_{\bf K})$,
where $e=id$. In addition $G^r_{u,\bf K}$ 
contains all mappings $f$ such that
$\phi \circ f(x)|_{U_j}=(id(x)+c'\bar Q_m(x)e_i)$ with $n=Ord(m)
\in \bf N$, $m \in \bf
N_o^n$, $i \in \bf N$, $c'\in \bf K$ and $\rho _a^t(f,id)\le r$ (that is,
for sufficiently small $|c'|_{\bf K}$
there is satisfied the following inequality 
$\| c'\bar Q_m\| _{C_{0,a}(t,U_j
\to X_{u,\bf K})}\le p^{-2}$). 
Therefore, the closure in $Vect_0(t,M)$ of the
$\bf K$-linear span of $\bigcup_uT_eG^r_{u,\bf K}$ 
coincides with $T_eG(t,M_{\bf K})$,
which follows from Kaplansky Theorem A.4 \cite{sch1}.
Evidently, $T_eG^r_{u,\bf K}$ is infinite-dimensional over $\bf K$. 
\par  Let us take $G^n_{u,\bf K}:=\{ f\in G:$ $supp(f)\subset U^{E_n},
f|_{U^{E_n}}\in C_{0,a}(t,U^{E_n}\to M_{u,\bf K})$, 
$\rho _a^t(f|_{U^{E_n}}, id) \le n \} $, 
where $U^E:=\bigcup_{j\in E}U_j,$ $E_n=(1,...,n)$, $n\in \bf N$, 
since $M_{u,\bf K}$ is separable. Each subgroup 
$G^n_{u,\bf K}$ is compact in $G$.
Since $\bigcup_{n\in \bf N} \{ f\in G:$ $supp(f)\subset U^{E_n} \} $
is dense in $G$, then $\bar G_a(t,M):=\bigcup_{n,u,\bf K}
G^n_{u,\bf K}$ is dense in $G$. If $f,$
$g\in \bar G_a(t,M)$, then there exists $n$ with $supp(f)\cup supp(g)
\subset U^{E_n}$ and $g^{-1}\circ f\in \bar G_a(t,M)$, hence
$\bar G_a(t,M)$ is the subgroup in $G$.
\par {\bf 3.7. Remarks.} If $M=B(X,0,1)$ for a normed space $X$, 
then $G_a(t,M_{u,\bf K})$ is a projective limit of discere groups 
$G_a(t,M_{u,\bf K})/B(G_a(t,M_{u,\bf K}),e,p^{-l})$
of polynomial bijective surjective mappings
${\tilde f}: k^n\to k^n$ of finite rings $k=B({\bf K},0,1)/
B({\bf K},0,p^{-l})$, since a series 
$f(x)=\sum_{m,i}a(m,f^i){\bar Q_m}(x)e_i$
in the Amice base ${\bar Q_m}e_i$
for each $f\in C^{id}_{0,a}(t,M_{u,\bf K}\to M_{u,\bf K})$ 
produces a finite sum 
${\tilde f}({\tilde x})=
\sum_{m,i}{\tilde a}(m,f^i){\tilde Q}_m({\tilde x})e_i , $
where $\pi _l: {\bf K}\to {\bf K}/B({\bf K},0,p^{-l})$
is the quotient mapping, $x\in B({\bf K},0,p^{-l})$,
${\tilde x}=\pi _l(x)$, $a(m,*)\in \bf K$, ${\tilde a}(m,*)=
\pi _l(a(m,*))$, ${\tilde Q}_m({\tilde x})=\pi _l
({\bar Q_m}(x))$, $l\in \bf N$ and $n\in \bf N$ 
depends on $l$, balls in $G_a(t,M_{u,\bf K})$
are given relative to the ultrametric $\rho _0^t$ in $G(t,M).$
If $X_{u,\bf K}$ is finite-dimensional over $\bf K$,
then $n$ are bounded by $dim_{\bf K}X_{u,\bf K}$.
\par Besides profinite subgroups
given in \S 3.6 there are classical subgroups over 
the non-Archimedean field $\bf F$ contained in
the diffeomorphism groups. In particular 
subgroups preserving vector fields are important
for quantum mechanics. Let $M$ be an analytic manifold, which is
is a clopen subset in
$B(X,0,r)$, where $r>0$, and $X=\bf F^n$. For a covector field
$\tilde A:=\{ A_{\alpha }(x):$ $\alpha =1,...,n \} $ on $M$
a differential 1-form 
\par $(i)$ $A=A_{\alpha }(x)dx^{\alpha }$
is called a potential structure,
where summation is by $\alpha =1,...,n$.  It is called analytic
if $\tilde A\in C(an_r,M\to {\bf F^n})$. It is called non-degenerate, if 
\par $(ii)$ $det(F_{\alpha ,\beta })\ne 0$ for each $x$, where 
\par $(iii)$ $F=dA=F_{\alpha ,\beta }dx^{\alpha }\Lambda dx^{\beta }/2$. 
We consider $g\in Diff(an_r,M)$ and $y^{\alpha }=g^{\alpha }(x)$, 
$x=(x^1,...,x^n)\in M$, $x^i\in \bf F$. If 
\par $(iv)$ $A_{\alpha }=A_{\mu }\partial g^{\mu }/\partial x^{\alpha }$ or 
\par $(v)$ $F_{\alpha ,
\beta }=F_{\mu ,\nu }(\partial g^{\mu }/\partial x^{\alpha })
(\partial g^{\nu }/\partial x^{\beta })$. The groups of such
$g$ are denoted by $G_A$ or $G_F$ respectively and are called 
a potential group or a symplectic
group respectively. Corresponding to them Lie algebras of vector
fields are denoted by $\sf L_A$ and $\sf L_F$.
There are accomplished analogs
of Proposals 1 and 2 \cite{parin}, since $G_A\subset G_F$ and
$\sf L_A\subset L_F$. Let $n=2m$, $m\in \bf N$ and 
\par $(vi)$ $A=c_{\alpha ,\nu }x^{\nu }dx^{\alpha }$, where $c_{\alpha ,
\nu }=-c_{\nu ,\alpha }=const $ and $det(c_{\alpha ,\nu })\ne 0$; or
\par $(vii)$ $A=A_{\alpha }dx^{\alpha }$, $A_{\alpha }=\sum _{k
=1}^{\infty }c^k_{\alpha ,\nu _1,...,\nu _k}x^{\nu _1}...x^{\nu _k}$,
where $c^k_{\alpha ,\nu _1,...,\nu _k}\in \bf F$, $c^1_{\alpha ,\nu }=
-c^1_{\nu ,\alpha }$ for each $\alpha ,\nu =1,...,N$,
$det(c^1_{\alpha ,\nu })\ne 0$. Then 
\par $(viii)$ $dim_{\bf F}{\sf L_A}
=n(n+1)/2$, $G_A=Sp(2m,{\bf F}):=\{ g\in GL(2m,{\sf F})|$
$g^t\epsilon g=\epsilon \} $ is the symplectic group,
where $g^t$ denotes the transposed matrix;
\par $(ix)$ $dim_{\bf F}{\sf L_A} \le n(n+1)/2 $.
\par To verify this let us consider at first $c_{\alpha ,
\mu }=\epsilon _{
\alpha ,\mu }$, where $\epsilon _{\alpha ,\alpha +1}=1$, $\epsilon
_{\alpha +1,\alpha }=-1$ for $\alpha =1,...,n-1$, $\epsilon _{
\alpha ,\beta }=0$ for others $(\alpha ,\beta )$. Therefore, 
there are true analogs of Formulas (10-13) \cite{parin} 
with $a^{\mu }_{\nu _1,
...,\nu _k}\in \bf F$. The matrices $B^{(k),\sigma }_{\alpha ,\nu ,
\mu }$ in Lemma 1 \S 2 have integer elements, consequently,
there are true analogs of Formulas (15,16) for the field $\bf F$,
since an analytic vector field $\xi $ is in $\sf L_A$ if and only if
$\xi ^{\mu }\partial _{\mu }A_{\alpha }+$ $A_{\mu }\partial _{
\alpha }\xi ^{\mu }=:L_{\xi }A_{\alpha }=0$.
In general the form $A$ can be reduced to $A=- \lambda
\epsilon _{\alpha ,\nu }x^{\nu }dx^{\alpha }/2$ by some operator
$j\in GL(n,{\bf F})$, where $\lambda \in \bf F$
and $j(B({\bf F^n},0,1))=B({\bf F^n},0,1)$. 
Theorem 2 \cite{parin} can also be modified,
but should be rather lengthy.
\par In $G(t,M)$ for $0\le t\le \infty $
there are also subgroups isomorphic with
the additive group $B(X_u,0,r)$, elements $f$ of which act as translations
of a subset $V$ of $M$ diffeomorphic with $B(X_u,0,r)$ 
and $f|_{M\setminus V}=id$. Using disjoint coverings 
of $M_{u,\bf K}$ by balls
we get, that $Diff(t,M)$ contains subgroups
isomorphic with symmetric groups $S_n$, where either $n\in \bf N$ 
or $n=\infty $ for non compact $M$.
Also $Diff(t,M)$ contains a subgroup diffeomorphic with  
$W:=\{ f:$ $f|_V\mbox{ has an extension}$ ${\tilde f}\in GL(X_u)$
$\| {\tilde f} -I \|_{X_u}<1,$ $f|_{M\setminus V}=id \} $,
where $GL(X_u)$ is the general linear group.
\section{ Decompositions of representations and induced representations.}
\par  {\bf 4.1.} Let $G=G(t,M)$ be defined as in 
\S 2.4 and \S 3.5 with $0\le t <\infty $ and $T: G\to U(H)$
be a strongly continuous unitary representation, 
where ${\bf Q_p}\subset {\bf F}\subset \bf C_p$,
$U(H)$ is a unitary group of a 
Hilbert space $H$ over $\bf C$. The unitary group is in the 
standard topology inherited from the space $L(H)$ of continuous 
linear operators $A: H\to H$ in the operator norm topology.
\par {\bf Theorem.}  {\it The representation $T$ 
up to the unitary equivalence $T\mapsto A^{-1}TA$ is decomposable 
into the direct 
integral $T_g=\int_JT_g(y)dv(y)$ of irreducible representations
$G\ni g\mapsto T_g(y)\in U(H_y)$, where $H_y$ are Hilbert subspaces of $H$,
$y\in J$, $v$ is a $\sigma $-additive measure on a compact 
Hausdorff space $J$, $A$ is a fixed unitary operator.}
\par {\bf Proof.} In view of Theorems 3.6 there exists a family of
compact subgroups $G^n_{u,\bf K}$ in $(G,V(G))$ for which 
$G^n_{u,\bf K}\subset G^{n+1}_{u,\bf K}$ 
for each $n$ and $N:={\bar G}_a(t,M)$
is dense in $(G,V(G))$, where $V(G)$ denotes the topology of $G$.
Then 
$$T_g|_{G^n_{u,\bf K}}=\int_{J(n,u,{\bf  K})}T_g(n,u,{\bf K};
y)v_{n,u,\bf K}(dy)$$
for each $n\in \bf N$, a pseudoultranorm $u$ in $X$ and a local subfield
${\bf K}\subset \bf F$, where $v_{n,u,\bf K}$ are measures on compact
spaces $J(n,u,{\bf K})$, $T_g(n,u,{\bf K};y)$ are finite-dimensional 
irreducible representations,
$y\in J(n,u,{\bf K})$, $g\in G^n_{u,\bf K}$. 
\par  There is the consistent family $T_g(n,u,{\bf K};y)$ such that
$v_{n,u,\bf K}$-almost everywhere in $J(n+1,u',{\bf K}')$ 
the restriction $T_g(n+1,u',{\bf K'};y)|_{G^n_{u,\bf K}}$
is a finite direct sum of $T_g(n,u,{\bf K};y)$ 
with the corresponding $y$, where $u(a)\le u'(a)$
for each $a\in X$, ${\bf K}\subset \bf K'$. Therefore, there are continuous
mappings $z(-n,u,{\bf K};-n',u',{\bf K'})$ from 
$J(n,u,{\bf K})$ into $J(n',u',{\bf K'})$ for each $n<n'$, $u\le u'$ 
and ${\bf K}\subset \bf K'$ such that
$v_{n',u',\bf K'}(Y)=v_{n,u,\bf K}(z^{-1}(-n,u,{\bf K};-n',u',{\bf K'})(Y))$ 
for each $Y$ in the Borel $\sigma $-filed $Bf(J(n',u',{\bf K'}))$,
where $v_{n,u,\bf K}$ are non-negative measures.
For each $ \xi , \eta \in  H$ with $ \| \xi \| = \| \eta \| =1$
we have $| (\xi ,T_gy)|\le 1$ and 
$$| \int_{J(n,u,{\bf  K})} (\xi ,T_g(n,u,{\bf K};y) \eta )
v_{n,u,\bf K}(dy)|\le 1.$$
Consequently, $T_g|_N=\int_J T_g(y)v(dy)$, where the projective
limit of compact spaces 
$J=pr-\lim \{ J(n,u,{\bf K}); z(-n,u,{\bf K};
-n',u'{\bf K'}); \{ (n,u,{\bf K}) \} \} $
is compact (see also \S 2.5 \cite{eng})
and the projective limit of measures $v=pr-lim \{ v_{n,u,\bf K} \}$
is the measure on $(J,Bf(J))$, and $T_g(y): N
\to U(H_y)$ are irreducible for $v$-almost every $y\in J$.
Therefore, 
$$T_g=\int_J T_g(y)v(dy) \mbox{ and } H=\int_J H_yv(dy),$$
where $T_g(y): G\to U(H_y)$ is an irreducible unitary representation
for $v$-almost each $y\in J$, $H_y$ are complex Hilbert subspaces of $H$ 
(see \cite{nai} and \S 22.8 \cite{hew}). 
\par {\bf 4.2.} Let $G:=G(t,M)$ be given by \S \S 2.4 and 3.5, where 
${\bf Q_p}\subset {\bf F}\subset \bf C_p$. Suppose
$W: G\to IS(H)$ is the regular representation (for $H$ over a local field 
${\bf L}\supset Q_s$, $s\ne p$) given by the formula $U_gf(x):=f(g^{-1}x)$,
where $H$ is a Banach space of uniformly continuous bounded
functions $f: G\to \bf L$ with a norm 
$\| f \| :=\sup_{x\in G} |f(x)|_{\bf L}$,
$IS(H)$ is a group of isometries of $H$ with a metric induced by 
an operator norm of continuous $\bf L$-linear operators $V$, $V: H\to H$. 
\par {\bf Theorem.} {\it There exists $A\in IS(H)$
such that $AWA^{-1}$ is decomposable into a direct sum
of irreducible representations $W_j$. Moreover,
each irreducible representation $T: G\to IS(E)$ for a Banach space 
$E$ over $\bf L$ is equivalent to some $W_j$.}
\par   {\bf Proof.} It may be directly verified that $W$ is strongly
continuous. This means that for each $c>0$ and $f\in H$ there is a
neinghbourhood $V\ni id$ such that
$\| W_gf-f\| \le c$ for each $g\in V$. Let the compact subgroups
$G^n_{u,\bf K}$ be the same
as in the proof of Theorem 4.1. They $s$-free, that is, 
for each clopen subgroup
$E'$ its index $[G^n_{u,\bf K}: E']$ is
not divisible by $s$ (\cite{roo,rosc}). 
It follows from the consideration of local
decompositions of elements in $G^n_{u,\bf K}$ by ${\bar Q}_me_i$
and from the fact that $B({\bf K},0,r)$ are the $s$-free additive
groups for each $0<r<\infty $, $n \in \bf N$.
In addition $E'$ contains an open compact
subgroup which is normal
in $G^n_{u,\bf K}$ due to Pontryagin lemma (see \S 8.1 \cite{roo}).
Therefore, on $G^n_{u,\bf K}$ the Haar measure exists with values in
$\bf Q_s$ due to Monna-Springer theorem 8.4\cite{roo}.
In view of Theorem 2.6 and Corollary 2.7 \cite{rosc} each
strongly continuous representation
$\tilde T: G^n_{u,\bf K}\to IS(H)$ is decomposable into the direct sum
of irreducible representations. On the other hand, ${\bar G}_a(t,M)$ 
is dense in $G.$ The last statement of this
theorem follows from the fact that for compact groups each
$T: G^n_{u,\bf K}\to IS(H_T)$
is equivalent to some irreducible component of the regular
representation, where $H_T$ is a Banach space over $\bf L$.
\par {\bf 4.3. Remark.} Let $\mu $ be a Borel regular Radon
non-negative quasi-invariant measure on a diffeomorphism group $G$ 
relative to a dense subgroup $G'$ with a continuous quasi-invariance factor
$\rho _{\mu }(x,y)$ on $G'\times G$ and
$0<\mu (G)<\infty $ \cite{lutmf99}.  
Suppose  that $V: S\to U(H_V)$ is a strongly continuous unitary 
representation of a closed subgroup $S$ in $G'$.
There exists a Hilbert space $L^2(G,\mu ,H_V)$ of equivalence classes 
of measurable functions $f: G\to H_V$ with a finite norm
$$(1)\mbox{ }\| f \| :=(\int_G \| f(g) \|^2_{H_V}\mu (dg))^{1/2}<\infty .$$
Then there exists a subspace $\Psi _0$ of functions
$f\in L^2(G,\mu ,H_V)$ such that $f(hy)=V_{h^{-1}}f(y)$ for each
$y\in G$ and $h\in S$, the closure of $\Psi _0$ in
$L^2(G,\mu ,H_V)$ is denoted by $\Psi ^{V,\mu }$. For each 
$f\in \Psi ^{V,\mu }$
there is defined a function
$$(2)\mbox{ }(T^{V,\mu }_xf)(y):=\rho _{\mu  }^{1/2}(x,y) f(x^{-1}y),$$
where $\rho _{\mu }(x,y):=\mu _x(dy)/\mu (dy)$ is a quasi-invariance 
factor for each $x\in G'$ and $y\in G$, $\mu _x(A):=
\mu (x^{-1}A)$ for each Borel subset $A$ in $G$.
Since $(T^{V,\mu }_xf)(hy)=V_{h^{-1}}((T_xf)(y))$, then $\Psi  ^{V,\mu }$
is a $T^{V,\mu }$-stable subspace. Therefore, $T^{V,\mu }: G'\to U(\Psi 
^{V,\mu })$
is a strongly continuous unitary representation,
which is called induced and denoted by $Ind_{S\uparrow G'}(V)$.
\par If $S=\lim \{ S_{\alpha }, \pi ^{\alpha }_{\beta }, \Omega \} $
is a profinite subgroup of $G$, for example, $G^n_{u,\bf K}$
(see \S \S 3.6, 3.7) and $V$ is irreducible, then
$H_V$ is finite-dimensional and $V^{-1}(I)=:W$ is a clopen normal 
subgroup in $S$, where $\pi ^{\alpha }_{\beta }: S_{\alpha }\to 
S_{\beta }$ are homomorphisms of finite groups $S_{\alpha }$ and 
$S_{\beta }$ for each $\alpha \le \beta $ in an ordered set
$\Omega $. Therefore, there exists $\alpha \in \Omega $ such that
$\pi _{\alpha }^{-1}(e_{\alpha })\subset W$, where $e_{\alpha }$
is the unit element in $S_{\alpha }$ and $\pi _{\alpha }:
S\to S_{\alpha }$ is a quotient homomorphism.
In view of Theorems 7.5-7.8 \cite{hew} there exists a representation
$V^{\alpha }: S_{\alpha }\to U(H_V)$ such that $V^{\alpha }\circ 
\pi _{\alpha }=V.$
\par {\bf 4.4.} Let $G$ be a diffeomorphism group
with a non-negative
quasi-invariant measure $\mu $ relative to a dense subgroup $G'$.
We can choose $G'$ such that each $G^n_{u,\bf K}$ is
a compact subgroup of $G'$. Suppose that there are two closed
subgroups $K$ and $N$ in $G$ such that $K':=K\cap G'$ and $N'=N\cap G'$
are dense subgroups in $K$ and $N$ respectively. We say that
$K$ and $N$ act regularly in $G$, if there exists a sequence 
$\{ Z_i:$ $i=0,1,... \} $ of Borel subsets $Z_i$ satisfying two conditions:
\par $(i)$ $\mu (Z_0)=0$, $Z_i(k,n)=Z_i$ for each pair
$(k,n)\in K\times N$ and each $i$;
\par $(ii)$ if an orbit $\sf O$ relative to the action of
$K\times N$ is not a subset of $Z_0$, then
${\sf O}=\bigcap_{Z_i\supset \sf O}Z_i$, where $g(k,n):=k^{-1}gn$.
Let $T^{V,\mu }$ be a representation of
$G'$ induced by a unitary representation $V$ of $K'$ and
a quasi-invariant measure $\mu $ as in \S 4.3. 
We denote by $T^{V,\mu }_{N'}$ a restriction of $T^{V,\mu }$ on $N'$
and by $\sf D$ a space $K\setminus G/N$ of double coset classes $KgN$.
\par {\bf Theorem.} {\it There are a unitary operator $A$ on 
$\Psi ^{V,\mu }$
and a measure $\nu $ on a space $\sf D$ such that 
$$(1)\mbox{ }A^{-1}T^{V,\mu }_nA=\int_{\sf D}T_n(\xi )d\nu (\xi )$$
for each $n\in N'$.
$(2).$ Each representation $N'\ni n\mapsto T_n(\xi )$ in 
the direct integral decomposition $(1)$ is defined relative to 
the equivalence of a double coset class $\xi $. For a subgroup
$N'\cap g^{-1}K'g$ its representations 
$\gamma \mapsto V_{g\gamma g^{-1}}$ are equivalent for each $g\in G'$
and representations $T_{N'}(\xi )$ can be taken up to their equivalence
as induced by $\gamma \mapsto V_{g\gamma g^{-1}}$.}
\par {\bf Proof.} A quotient mapping 
$\pi _{\sf X}: G\to G/K=:\sf X$
induces a measure $\hat \mu $ on $\sf X$ such that
${\hat \mu }(E)=\mu (\pi ^{-1}_{\sf X}(E))=:(\pi _{\sf X}^*\mu )(E)$ 
for each Borel subset 
$E$ in $\sf X$. In view of Radon-Nikodym theorem II.7.8 \cite{fell}
for each $\xi \in \sf D$
there exists a measure $\mu _{\xi }$ on $\sf X$ such that
$$(3)\mbox{ }d{\hat \mu }(x)= d\nu (\xi )d\mu _{\xi }(x),$$ 
where $x\in \sf X$, $\nu (E):=(s^*\mu )(E)$
for each Borel subset $E$ in $\sf D$, $s: G\to \sf D$ is
a quotient mapping.
\par For each $m\in \bf N$, a pseudoultranorm $u$ in $X$ and
a local subfield $\bf K$ in $\bf F$ a subgroup
$G^m_{u,\bf K}$ is compact in $G$, hence there exists
a topological retraction $r_{m,u,\bf K}: G\to G^m_{u,\bf K}$ 
(that is, $r_{m,u,\bf K}\circ r_{m,u,\bf K}=r_{m,u,\bf K}$ and
$r_{m,u,\bf K}$ is continuous 
and $r_{m,u,\bf K}|_{G^m_{u,\bf K}}=id$).
This retraction induces a measure $\mu _{m,u,\bf K}=
r_{m,u,\bf K}^*\mu $ on $G^m_{u,\bf K}$.
It is equivalent to a Haar measure on $G^m_{u,\bf K}$, since
it is quasi-invariant relative to $G^m_{u,\bf K}$ (see \S
VII.1.9 in \cite{boui}). In view of \S 26 \cite{nari} and Formula $(3)$
the Hilbert space $H^V:=L^2({\sf X},{\hat \mu },H_V)$ has a decomposition 
into a direct integral 
$$(4)\mbox{ } H^V=\int_{\sf D}H(\xi )d\nu (\xi ),$$
where $H_V$ denotes a complex Hilbert space of the representation
$V: K'\to U(H_V)$. Therefore, 
$$ \| f \| ^2_{H^V}=\int_{\sf D}\| f \| ^2_{H(\xi )}d\nu (\xi ).$$
In view of \S 32.2 from Chapter VI \cite{nai}
each irreducible representation of a compact group $Y$
can be realized as a representaion in some minimal left ideal
of a ring $L^2(Y, \lambda ,{\bf C})$, where $\lambda $ is a 
Haar measure on $Y$. 
From Formulas $(4)$ and $4.3.(1,2)$ we get the first statement 
of this theorem for a subspace $\Psi ^{V,\mu }$ of $H^V$.
\par If $f\in L^2({\sf X},{\hat \mu },H_V)$, then
$\pi ^*_{\sf X}f:=f\circ \pi _{\sf X}\in L^2(G,\mu ,H_V)$.
This induces an embedding $\pi ^*_{\sf X}$ of $H^V$ into
$\Psi ^{V,\mu }$. Let $\sf F$ be a filterbase of neighbourhoods $A$ 
of $K$ in $G$ such that $A=\pi ^{-1}_{\sf X}(S)$, where
$S$ is open in $\sf X$, hence $0<\mu (A)\le \mu (G)$
due to quasi-invaraince of $\mu $ on $G$ relative to $G'$.
Let $\psi \in \xi \in \sf D$, then $\psi =Kg_{\xi }$, where
$g_{\xi }\in G$, hence $\psi =\psi (N\cap g_{\xi }^{-1}Kg_{\xi })$.
In view of Formula $(3)$ for each $x\in N'$ and $\eta
=Kx$ we get
$\rho _{\mu _{\xi }}^{1/2}(\eta ,\xi )=
\lim_{\sf F}[\int_A \rho ^{1/2}(x,zg_{\xi })d\mu (z)/\mu (A)]$
(see also \S 1.6 \cite{eng}).
Therefore, $(a,T_x(\xi )b)_{H^V}=$ $\lim_{\sf F} [\int_A (\pi ^*a,$ 
$ \rho _{\mu }^{1/2}(x,zg_{\xi })(\pi ^*b)_x^{zg_{\xi }})_{
\Psi ^{V,\mu }}d\mu (z)/\mu (A)]$
for each $x\in N'$ and $a,b \in H^V$, where $f_z^h(\zeta ):=
f(z^{-1}h\zeta )$ for a function $f$  on $G$ and $h, z, \zeta \in G$. 
In view of the cocycle condition $\rho _{\mu }(yx,z)=\rho _{\mu }
(x,y^{-1}z)\rho _{\mu }(y,z)$ for each $x, y \in G'$ and $z\in G$ we get 
$T_{yx}(\xi )=T_y(\xi )T_x(\xi )$ for each $x, y\in N'$ and $T_x(\xi )$ 
are unitary representations of $N'$. Then $(a,T_{yx}(\xi )b)_{H^V}=$
$\lim_{\sf F} [\int_A (\pi ^*a,$ $V_{g_{\xi }yg_{\xi }^{-1}}[$
$\rho _{\mu }^{1/2}(x,zg_{\xi })(\pi ^*b)_x^{zg_{\xi }}])_{\Psi ^{V,\mu }}
d\mu (z)/\mu (A)]$ for each $y\in N'\cap g_{\xi }^{-1}K'g_{\xi }$.
Hence the representation $T_x(\xi )$ in the Hilbert space
$H(\xi )$ is induced by a representation 
$(N'\cap g_{\xi }^{-1}K'g_{\xi })\ni y\mapsto 
V_{g_{\xi }yg_{\xi }^{-1}}.$
\par {\bf 4.5.} Let $V$ and $W$ be two unitary representations of 
$K'$ and $N'$ (see \S 4.4). In addition let $K$ and $N$ be regularly 
related in $G$ and $V{\hat \otimes }W$ denotes an external tensor product of
representations for a direct product group $K\times N$.
In view of \S 4.3 a representation $T^{V,\mu }{\hat \otimes }T^{W,\mu }$
of an external product group ${\sf G}:=G\times G$ 
is equvalent with an induced representation
$T^{V{\hat \otimes }W,\mu \otimes \mu }$, where $\mu \otimes \mu $
is a product measure on $\sf G$.
A restriction of $T^{V{\hat \otimes }W,\mu \otimes \mu }$
on ${\tilde G}:=\{ (g,g):$ $g\in G \}$
is equivalent with an internal tensor product
$T^{V,\mu }\otimes T^{W,\mu }$.
\par {\bf Theorem.} {\it There exists a unitary operator $A$ 
on $\Psi ^{V{\hat \otimes }W,\mu \otimes \mu }$ such that
$$(1)\mbox{ }A^{-1}T^{V,\mu }\otimes T^{W,\mu }A=\int_{\sf D}
T(\xi )d\nu (\xi ),$$
where $\nu $ is an admissible measure on a space ${\sf D}:=
N\setminus G/K$ of double cosets.
\par $(2).$ Each representation $G'\ni g\mapsto T_g(\xi )$
in Formula $(1)$ is defined up to the equivalence of $\xi $ in $\sf D$.
Moreover, $T(\xi )$ is unitarily equivalent with
$T^{{\tilde V}\otimes {\tilde W},\mu \otimes \mu }$, 
where $\tilde V$ and $\tilde W$
are restrictions of the corresponding representations 
$y\mapsto V_{gyg^{-1}}$ and $y\mapsto W_{\gamma y\gamma ^{-1}}$ 
on $g^{-1}K'g\cap \gamma ^{-1}N'\gamma $,
$g, \gamma \in G'$, $g\gamma ^{-1}\in \xi $.}
\par {\bf Proof.} In view of \S 18.2 \cite{barut}
$P\setminus {\sf G}/{\tilde G}$ and $K\setminus G/N$ 
are isomorphic Borel spaces,
where $P=K\times N$. In view of Theorem 4.4
there exists a unitary operator $A$ on a subspace $\Psi ^{V
\hat \otimes W,\mu \otimes \mu }$ of the Hilbert space
$L^2({\sf G},\mu \otimes \mu ,H_V\otimes H_W)$ such that 
$$A^{-1}T^{V\hat \otimes W,\mu \otimes \mu }|_{\tilde G}A=
\int_{\sf D}T_{\tilde G}(\xi )d\nu (\xi ),$$
where each $T_{\tilde G}(\xi )$ is induced by
a representation $(y,y)\mapsto (V\hat \otimes W)_{(
g,\gamma )(y,y)(g,\gamma )^{-1}}$ of a subgroup
${\tilde G}'\cap (g,\gamma )^{-1}(K\times N)(g,\gamma )$,
the latter group is isomorphic with
$S:=g^{-1}K'g\cap \gamma ^{-1}N'\gamma ,$ that gives a representation
${\tilde V}{\hat \otimes }{\tilde W}$ of a subgroup $S\times S$ in $\sf G$.
Therefore, we get a representation $T^{{\tilde V}{\hat \otimes }{\hat W},
\mu \otimes \mu }$ equivalent with $Ind_{(S\times S)\uparrow {\sf G}'}
({\tilde V}{\hat \otimes }{\tilde W})|_{{\tilde G}'}$.

\end{document}